\documentclass[12pt,reqno]{amsart}
\usepackage{amsmath}
\usepackage{amsxtra}
\usepackage{amscd}
\usepackage{amsthm}
\usepackage{amsfonts}
\usepackage{amssymb}
\usepackage{eucal}

\newtheorem{theorem}{Theorem}[section]

\newtheorem{cor}[theorem]{Corollary}
\newtheorem{lem}[theorem]{Lemma}
\newtheorem{prop}[theorem]{Proposition}

\theoremstyle{definition}
\newtheorem{defn}[theorem]{Definition}

\theoremstyle{remark}
\newtheorem{rem}[theorem]{Remark}
\theoremstyle{remark}
\newtheorem{example}[theorem]{Example}

\numberwithin{equation}{section}

\newcommand{\nc}{\newcommand}

\nc{\prp}{\perp\hskip -.1cm\perp}

\nc{\on}{\operatorname}
\nc{\ch}{\mbox{ch}}
\nc{\Z}{{\mathbb Z}}
\nc{\C}{{\mathbb C}}
\nc{\Pb}{{\mathbb P}}
\nc{\pone}{{\mathbb C}{\mathbb P}^1}
\nc{\pa}{\partial}
\nc{\F}{{\mathcal F}}
\nc{\arr}{\rightarrow}
\nc{\larr}{\longrightarrow}
\nc{\al}{\alpha}
\nc{\ri}{\rangle}
\nc{\lef}{\langle}
\nc{\W}{{\mathcal W}}
\nc{\la}{\lambda}
\nc{\ep}{\epsilon}

\nc{\su}{\widehat{{\mathfrak sl}}_2}
\nc{\sw}{{\mathfrak s}{\mathfrak l}}

\nc{\g}{{\mathfrak g}}
\nc{\h}{{\mathfrak h}}
\nc{\n}{{\mathfrak n}}
\nc{\N}{\widehat{\n}}
\nc{\G}{\widehat{\g}}
\nc{\De}{\Delta_+}
\nc{\gt}{\widetilde{\g}}
\nc{\Ga}{\Gamma}
\nc{\one}{{\mathbf 1}}
\nc{\z}{{\mathfrak Z}}
\nc{\zz}{{\mathcal Z}}
\nc{\Hh}{{\mathcal H}_\beta}
\nc{\qp}{q^{\frac{k}{2}}}
\nc{\qm}{q^{-\frac{k}{2}}}
\nc{\La}{\Lambda}
\nc{\wt}{\widetilde}
\nc{\qn}{\frac{[m]_q^2}{[2m]_q}}
\nc{\cri}{_{\on{cr}}}
\nc{\kk}{h^\vee}
\nc{\sun}{\widehat{\sw}_N}
\nc{\hh}{\widehat{\mathfrak h}}
\nc{\HH}{{\mathcal H}_{q,t}}
\nc{\ca}{\wt{{\mathcal A}}_{h,k}(\sw_2)}
\nc{\gl}{\widehat{{\mathfrak g}{\mathfrak l}}_2}
\nc{\el}{\ell}
\nc{\s}{{\mathbf s}}
\nc{\bi}{\bibitem}
\nc{\om}{\omega}
\nc{\WW}{\W_\beta}
\nc{\scr}{{\mathbf S}}
\nc{\ab}{{\mathbf a}}
\nc{\rr}{r}
\nc{\ol}{\overline}
\nc{\con}{qt^{-1} + q^{-1}t}
\nc{\den}{q^{\el-1} t^{-\el+1}+ q^{-\el+1} t^{\el-1}}
\nc{\ds}{\displaystyle}
\nc{\B}{B}
\nc{\A}{{\mathbb A}}
\nc{\GG}{{\mathcal G}}
\nc{\UU}{{\mathcal U}}
\nc{\MM}{{\mathcal M}}
\nc{\CC}{{\mathcal C}}
\nc{\GL}{{}^L G}
\nc{\dzz}{\frac{dz}{z}}
\nc{\Res}{\on{Res}}
\nc{\rep}{{\mathcal R}ep \;}
\nc{\uqg}{U_q \G}
\nc{\uqgg}{U_q \g}
\nc{\Fq}{{\mathbb F}_q}

\nc{\stimes}{\ltimes}
\nc{\K}{\hat{\mathcal K}}
\nc{\Ql}{\ol{\mathbb Q}_\ell}

\nc{\ga}{\gamma}
\nc{\PL}{{}^L P}
\nc{\E}{\mc E}
\nc{\mc}{\mathcal}
\nc{\mbf}{\mathbf}
\nc{\bb}{{\mathfrak b}}
\nc{\OO}{{\mc O}}
\nc{\Po}{{\mc P}}
\nc{\V}{{\mc V}}
\nc{\yy}{{\mc Y}}
\nc{\M}{\mathcal M}
\nc{\Coh}{{{\mathcal C}oh}}
\nc{\Cohn}{\Coh_n}
\nc{\f}{{\mathcal F}}
\nc{\si}{_E}
\nc{\Gaf}{{\mathbb G}_{a,\Fq}}
\nc{\KK}{{\mathfrak k}}

\nc{\PCr}{{ \bs P  (\C[x])^r   }}
\nc{\PCN}{{ \bs P  (\C[x])^N   }}
 
\nc{\sN}{sl_{2N+1}}

\nc{\Pzr}{{ \bs P(\C((x-z)))^r}}

\nc{\PzN}{{ \bs P(\C((x-z)))^N}}

\newcommand{\bean}{\begin{eqnarray}}
\newcommand{\eean}{\end{eqnarray}}
\newcommand{\be}{\begin{displaymath}}
\newcommand{\ee}{\end{displaymath}}
\newcommand{\bea}{\begin{eqnarray*}}   
\newcommand{\eea}{\end{eqnarray*}}
\newcommand{\bs}{\boldsymbol}
\newcommand{\Ref}[1]{{$($\ref{#1}$)$}}

\begin{document}
\title[Self-self-dual spaces of polynomials]
{Self-self-dual spaces of polynomials}
\author[Lev Borisov and Evgeny Mukhin]
{Lev Borisov and Evgeny Mukhin}
\thanks{Research of L.B. is supported in part by NSF grant DMS-0140172.
Research of E.M. is supported in part by NSF grant DMS-0140460.}
\address{L.B.:  Mathematics Department, University of Wisconsin-Madison,
480 Lincoln Dr, Madison, WI 53706-1388, USA, \newline borisov@math.wisc.edu}
\address{E.M.: Department of Mathematical Sciences, Indiana University -
Purdue University Indianapolis, 402 North Blackford St, Indianapolis,
IN 46202-3216, USA, \newline mukhin@math.iupui.edu}

\begin{abstract}
A space of polynomials $V$ of dimension 7 
is called self-dual if the divided Wronskian of
any 6-subspace is in $V$. A self-dual space $V$ has a natural inner
product. The divided Wronskian of any isotropic 3-subspace of $V$ 
is a square of a polynomial. We call $V$ self-self-dual if the square
root of the divided Wronskian of any isotropic 3-subspace is again
in $V$. We show that the self-self-dual spaces have a natural
non-degenerate skew-symmetric 3-form defined in terms of Wronskians.

We show that the self-self-dual spaces correspond to $G_2$-populations
related to the Bethe Ansatz of the Gaudin model of type $G_2$ 
and prove that a $G_2$-population is isomorphic to the $G_2$ flag variety.
\end{abstract}

\bigskip

\maketitle

\section{Introduction}
The Bethe equation is the main equation in the Bethe Ansatz method of
diagonalizing the Hamiltonians of many integrable systems of
mathematical physics. Given a solution of the Bethe equation, one
explicitly constructs an eigenvector of the Hamiltonian. 
This paper is related to the case of Gaudin model associated to $G_2$.
In this case, the Bethe equation is the system of algebraic equations
on complex variables $t_i^{(j)}$ with parameters
$z_i\in\C,m_s^{(j)}\in \Z_{>0}$: 
\bea\label{Bethe eqn}
-\sum_{s=1}^n\frac{m_s^{(1)}}{t_i^{(1)}-z_s}-
\sum_{k=1}^{l_2}\frac3{t_i^{(1)}-t_k^{(2)}}+\sum_{k=1,\ k\neq
  i}^{l_1}\frac2{t_i^{(1)}-t_k^{(1)}}=0,\qquad i=1,\dots,l_1,\notag\\
-\sum_{s=1}^n\frac{3m_s^{(2)}}{t_j^{(2)}-z_s}-
\sum_{k=1}^{l_1}\frac3{t_j^{(2)}-t_k^{(1)}}+\sum_{k=1,\ k\neq
  j}^{l_2}\frac6{t_j^{(2)}-t_k^{(2)}}=0, \qquad j=1,\dots,l_2.
\eea

A pair of polynomials $(y_1,y_2)\in(\Pb \C[x])^2$, is called
fertile with respect to polynomials $T_1,T_2$ if there exist
polynomials $\tilde y_1,\tilde y_2$ such that we have explicit Wronskians:
\be
W(y_1,\tilde y_1)=T_1y_2,\qquad W(y_2,\tilde y_2)=T_2y_1^3.
\ee 
The pair $(y_1,y_2)$ is called generic if $y_i(x)$
have no multiple roots and no common roots.

It is shown in \cite{MV} that 
zeroes of the pair of polynomials
$(\prod_{i=1}^{l_1}(x-t_i^{(1)}),\prod_{j=1}^{l_2}(x-t_j^{(2)}))$
satisfy the Bethe 
equation if and only if the pair is generic and fertile with respect
to $T_j(x)=\prod_{i=1}^n(x-z_s)^{m_s^{(j)}}$. We call such a pair a
Bethe pair.

Moreover, it is shown in \cite{MV} that if the pair $(y_1,y_2)$ is
a Bethe pair, then the pairs $(\tilde y_1,y_2)$ and $(y_1,\tilde
y_2)$ are fertile and therefore for almost all choices of $\tilde
y_1,\tilde y_2$ these pairs are Bethe pairs. Thus given one Bethe pair
we obtain a family of new Bethe pairs which in turn produce new Bethe
pairs, etc. The Zariski closure of 
all Bethe pairs obtained from $(y_1,y_2)$ is called a $G_2$-population
originated at $(y_1,y_2)$. 

It is conjectured in \cite{MV} that the number of $G_2$-populations
for generic $z_i$ equals the multiplicity of $L_{\La_{\infty}}$ in
$L_{\La_1}\otimes\dots\otimes L_{\La_n}$. Here $L_\La$ denotes an irreducible
finite-dimensional $G_2$-module with highest weight $\La$, 
$(\La_s,\al_i^\vee)=m_s^{(i)}$ and
$\La_\infty$ is the unique dominant weight in the Weyl group orbit of the
$G_2$-weight $\sum_{i=1}^n\La_i-l_1\al_1-l_2\al_2$ where $\al_i$ are
simple roots of a $G_2$ root system.
   
The $G_2$-populations are the original motivation of this paper. It
turns out that that the $G_2$-populations are in one-to-one correspondence
with special 7-dimensional spaces of polynomials, which we call
self-self-dual. 

\medskip

Let $V$ be a 7-dimensional space of polynomials. Let $U_i$ be the
greatest common monic divisor of Wronskians 
$\{W(v_1,\dots,v_i), v_j\in V\}$. 
We always assume that $U_1=1$. 
Then there exist unique polynomials $T_1,\dots,T_6$ such
that $U_i=T_{i-1}T_{i-2}^2\dots T_1^{i-1}$.

The polynomial $W^\dagger(v_1,\dots,v_i)=W(v_1,\dots,v_i)/U_i$ is called 
the divided Wronskian. 

The space $V$ is called self-dual if the space of all divided
6-Wronskians coincides with $V$. Self-dual spaces were studied in
\cite{MV}. Note that if $V$ is self-dual then $T_i=T_{7-i}$.
A self-dual space has a natural
non-degenerate symmetric bilinear form given by
\be
B(u,v)=W^\dagger(u,v_1,\dots,v_6), \qquad {\rm if}\ v=W^\dagger(v_1,\dots,v_6).
\ee
Moreover, any divided 3-Wronskian of an isotropic 3-space is a perfect
square. 

We call a self-dual space {\it self-self-dual} if 
\be
\{W^\dagger(v_1,v_2,v_3), v_i\in V, B(v_i,v_j)=0\}\ =\ \{v^2,v\in V\}.
\ee 
Note, that if $V$ is self-self-dual then $T_3=T_1$, so all $T_i$ are
expressed in terms of $T_1$ and $T_2$.
We prove that a self-self-dual space has a skew-symmetric
non-degenerate 3-form $w$ uniquely determined by the condition:
\be
w(v_1,v_2,v_3)=B(v,v),\qquad {\rm if}\
W^\dagger(v_1,v_2,v_3)=v^2\ {\rm and}\ B(v_i,v_j)=0.  
\ee

We identify the group $G_2$ with the subgroup of
orthogonal group $SO(V)$ which preserves $w$. It follows that 
the $G_2$ flag variety $G_2/B$ is identified
with the variety of $G_2$-isotropic flags 
$F=\{F_1\subset F_2\subset\dots \subset
F_7=V\}$ which have the properties 
$B(F_i,F_{7-i})=0$ and $F_3=\on{Ker}( w(F_1,\cdot,\cdot))$.  
We supply proofs of these and some other general facts on $G_2$ which
we failed to find in the literature.

Then we show that the first coordinates of a $G_2$-population span a
self-self-dual space $V$ and the $G_2$-population 
is isomorphic to the variety of $G_2$-isotropic flags in $V$. 
The isomorphism maps a flag $F$ to the pair $(F_1,W^{\dagger}(F_2))$. 

Let us say a few words about our methods. A self-dual space $V$ is
naturally a vector representation of $SO(7)$. However, we have an
explicit expression for the value of the 3-form only on 
isotropic triples of vectors in $V$. The lack of general formula makes it
difficult to prove the well-definedness of the 3-form and to compute
it. To overcome this problem, 
we consider the 8-dimensional spin representation $\hat V$ of
$Spin(7)$. Vectors in $V$ naturally act on $\hat V$. Moreover we
show that the set of isotropic triples in $V$ embeds in the
projectivization of $\hat V$ as a non-degenerate conic. This conic
defines an inner product $\hat B$ on $\hat V$ and the vectors of $V$ act by
skew-symmetric operators. 
It turns out that the 3-form can be defined
for all triples, isotropic or not, by the same formula $w(a,b,c)=\hat
B(abc\cdot p,p)$ where $p\in \hat V$ is also computed explicitly. 
In addition, the 2-form $B$ is given by the formula 
$B(v,w)=2\hat B(v\cdot p,w\cdot p)/\hat B(p,p)$. 
Those key observations allow us to finish all the proofs.

\medskip

The paper is constructed as follows. In Section \ref{sp emb} we 
prepare the main facts about spinor embedding of isotropic
Grassmannian of the vector representation of $SO(7)$. In Section \ref{G2
  sec} we collect and prove some general facts about $G_2$. In Section
\ref{ssd spaces} we define and study self-self-dual spaces of
polynomials. In Section \ref{populations} we show that
$G_2$-populations are in one-to-one correspondence with self-self-dual
spaces. In Section \ref{bases} we construct a special basis in a
self-self-dual space so that all divided 3-Wronskians, and the 3-form
have a canonical form.

\medskip
{\bf Acknowledgments.} E.M. is grateful to A. Varchenko for
stimulating discussions.

\medskip

\section{Spinor embedding}\label{sp emb}
In this section we define two important technical tools of this
paper. First we define and study the
spinor embedding of isotropic Grassmannian in the vector representation of
$SO(7)$ to the projectivization of the spin representation. Then we
proceed to the special maps from the spin representation to the vector
representation which we call invariant surjections.
\subsection{Spin representation}
We recall the basic statements about the Clifford algebra and
the algebraic group $Spin(7)$. We generally follow the book of Adams,
\cite{A}. However, our exposition is adapted to the problem
of embedding of the Grassmannian of isotropic 3-spaces 
in a seven-dimensional space with a symmetric non-degenerate form
into the projectivization of the spin representation of $Spin(7)$.

Let $V$ be a seven-dimensional complex vector space equipped with 
a non-degenerate quadratic form $Q$.  Let $B(\cdot , \cdot )$ 
be the bilinear form given by $B(v,w)=Q(v+w)-Q(v)-Q(w)$. Note that
$B(v,v)=2Q(v)$. 

By definition, the Clifford algebra $Cl(V)$ of
$V$ is the quotient of the free tensor algebra $\oplus_{n\geq 0}V^{\otimes n}$
by the two-sided ideal generated by $v\otimes v \oplus Q(v)$ for all $v\in V$.
It contains the even Clifford subalgebra $Cl^+(V)$ given by the the image
of the even graded components of the free tensor algebra. Dimensions
of $Cl(V)$ and $Cl^+(V)$ are $2^7$ and $2^6$ respectively.

The spinor group $Spin(V)$ is defined as the subgroup of the invertible
elements $g$ of $Cl^+(V)$ such that  $gVg^{-1}=V$.
The group  $Spin(V)$ is a double cover of $SO(V)\cong SO(7)$  with the map 
given by the above conjugation action on $V$, see \cite{A}.

To describe the irreducible representations of $Cl(V)$ and $Cl^+(V)$ 
we  identify these algebras with quotients of 
group algebras of finite $2$-groups as follows.
For any orthonormal basis $v_1,\dots,v_7$ 
of $V$ the images of the elements $v_i$ in $Cl(V)$
generate a group $H$ of order $2^8$. The elements of $H$ are of the form
$\pm \prod_{i=1}^7 v_i^{r_i}$ with $r_i\in \{0,1\}$. The defining relations
are $v_i^2 =-1$ and $v_iv_j=-v_jv_i$ for $i\neq j$. Let $H^+$
be the subgroup of $H$ given by the condition
$\sum_{i=1}^7 r_i=0\ ({\rm mod }\ 2)$. The group $H^+$ has index $2$ in $H$,
and the order of $H^+$ is $2^7$. We denote by the same letter $w$ the
central element $-1$ in $H$ and $H^+$, we have $w^2=1$.

The algebra 
$Cl(V)$  (resp. $Cl^+(V)$)
is the quotient of the group algebra of $H$ ($H^+$) by $w+1$.  
Consequently, complex representations of $Cl(V)$ (resp. $Cl^+(V)$)
are in one-to-one correspondence with 
representations of $H$ (resp. $H^+$) for which $w$ acts by $(-1)$. 

The quotient groups with respect to $\Z/2\Z$ subgroups generated by
 $w$, $H/\{1,w\}$ and $H^+/\{1,w\}$, are abelian $2$-groups 
of orders $2^7$ and $2^6$ respectively. Therefore, the group algebra
of $H$ (resp. $H^+$) 
has $2^7$ (resp. $2^6$) one-dimensional representations where $w$
acts by 1. Recall that the irreducible representations of a group
 algebra are in one-to-one 
correspondence with conjugacy classes in the group.

Since $7$ is an odd number, 
conjugacy classes in $H^+$ are $\{1\}$, $\{w\}$ $\{\pm h\}$, where 
$h\in H^+/\{1,w\}$, $h\neq 1$. 
Thus, there is only one additional representation
of $H^+$. Since the sum of squares of dimensions of all irreducible
representations equals the order of the group, the 
dimension of this additional representation $d$ 
satisfies $d^2+2^6=2^7$ which yields $d=8$. 
This produces a representation of $Cl^+(V)$ and $Spin(V)$ which
is called \emph{spin representation}. We denote this representation
$\hat V$.

The conjugacy classes of the
group $H$ are $\{1\}$, $\{w\}$, $\{w_1:=v_1v_2\ldots v_7\}$,
$\{w_2:=w w_1\}$, 
$\{\pm h\}$ where $h\in H/\{1,w\}$, $h\neq 1$, $h\neq \bar w_1$.
As a result, there are two additional representations
of $H$. On the other hand, we have two surjective group homomorphisms
$H\to H^{+}$ given by modding out the central element 
$w_1$ or $w_2$. These homomorphisms produce two non-isomorphic irreducible
representations of $H$ of dimension $8$. 

The group $H$ admits an outer automorphism $\alpha$
which sends $\prod_{i=1}^7 v_i^{r_i}$ to $w^{\sum_ir_i}
\prod_{i=1}^7 v_i^{r_i}$. This automorphism permutes the two
representations of $H$. Hence, the irreducible
$8$-dimensional spin representation $\hat V$ of $Cl^+(V)$ and $Spin(V)$ can
be lifted to a representation of the full Clifford algebra $Cl(V)$ in two
ways that differ by the automorphism $\alpha$ of $Cl(V)$ that
preserves $Cl^+(V)$ and multiplies odd degree elements of $Cl(V)$ by
$(-1)$. In particular, there is an action of elements of $V$ on the 
space $\hat V$, defined uniquely up to an overall sign. 

By the construction, the action of $Spin(V)$-module $V$ on $Spin(V)$-module 
$\hat V$ has the following equivariance property 
\bean\label{equiv}
g(v\cdot x)=(gv)\cdot(gx),
\eean 
where $g\in Spin(V)$, $v\in V$ and $x\in \hat V$. 

Now we describe the representation $\hat V$ explicitly.
Let $\hat V$  be the space of polynomials in the
odd variables $\hat v_5,\hat v_6,\hat v_7$. The space $\hat V$
has the basis  $\{\hat 1,\hat v_5,\hat v_6,\hat 
v_7,\hat v_5\hat v_6,\hat v_6\hat v_7,\hat v_5\hat v_7,\hat
v_5\hat v_6\hat v_7\}$.

Choose a  basis $v_1,\dots,v_7$ of $V$ with the property
$B(v_i,v_j)=(-1)^{i+1} \delta_{i+j}^8$. We will use 
the convention $Q(v)=\frac 12 B(v,v)$.
The algebra $Cl(V)$ is the associative $\mathbb C$-algebra
with generators $v_i$ and relations 
\bean\label{spinvone}
v_iv_j+v_jv_i=(-1)^{i}\delta_{i+j}^8.
\eean
Define an action of $Cl(V)$ in $\hat V$ as follows.
We let $v_5,v_6,v_7$ act by multiplications by $\hat v_5,\hat
v_6, \hat v_7$ respectively, $v_1$, $v_2$, $v_3$
by differentiations $\partial/\partial v_7$, $-\partial/\partial
v_6$,$\partial/\partial v_5$ 
respectively, and $v_4$ by $(-1)^{\rm deg}$,
where ${\rm deg}$ is the degree of the odd polynomial.

\begin{lem}\label{spin model}
The space $\hat V$ is an 8-dimensional irreducible representation
of $Cl(V)$. In particular, $\hat V$ is the spin representation of $Spin(V)$. 
\end{lem}
\begin{proof}
It is an easy check that equations \eqref{spinvone} are satisfied. The
irreducibility is obvious.
\end{proof}

We also need the decomposition of the second symmetric power
of $\hat V$ into irreducible representations of $Spin(V)$, see \cite{A}.
\begin{prop}\label{symtwo} We have the isomorphism of
  $Spin(V)$-modules 
defined uniquely up to scalars
\be
Sym^2(\hat V)\cong \C \oplus \Lambda^3(V),
\ee
where $\C$ is the trivial representation.
\end{prop}
We also observe that while $SO(V)$ acts on $\hat V$ only projectively,  
the $Sym^2(\hat V)$ is an honest representation of $SO(V)$, and
therefore the above proposition is also a decomposition of 
representations of the group $SO(V)$. 

\medskip


\subsection{Spinor embedding of isotropic Grassmannian in $\Pb \hat V$}

A subspace $U\subset V$ is called {\it isotropic} if the the quadratic
form $Q$ vanishes on $V$.  Let ${IG}(3,V)$ 
denote the Grassmannian of isotropic three-dimensional  subspaces
in $V$.

Let $U\subset V$ be an isotropic subspace of $V$ of dimension 3. 
Let $L_{U}\subset \hat V$ be the common kernel of elements of $U$
acting on $\hat V$: 
\be 
L_{U}:=\{x\in \hat V\ | \ v\cdot x=0,\ v\in U\}.
\ee
Here $v\cdot x$ is the action of $v\in V$ on $x\in\hat V$.
The action is defined up to a sign, therefore the space $L_{U}$ is
well defined.

For any complex vector space $W$ we denote by $\Pb W$ the
corresponding projective space, thought of as the space of lines in
$W$. For any non-zero $w\in W$ we have $\C w\in \Pb W$. 

\begin{theorem}\label{rhomap}
For every isotropic 3-subspace $U\subset V$ the space
$L_{U}$ is one-dimensional. The map
$\rho:\ {IG}(3,V)\to \Pb\hat V$, sending $U\mapsto L_U$ is a
$Spin(V)$-equivariant embedding whose image is 
a smooth degree two hypersurface in $ \Pb\hat V$.
\end{theorem}
We will call the map $\rho$ the \emph{spinor embedding}.
\begin{proof}
Let $v_1,v_2,v_3$ be a basis of $U$. We extend it to a basis of $V$
such that $B(v_i,v_j)=(-1)^{i+1}\delta_{i+j}^8$ and use the model of $\hat
V$ as in Lemma \ref{spin model}.

It is an easy calculation to see that $L_U=\mathbb C\hat 1$ 
and is therefore one-dimensional and the first part of theorem is proved.

Now that the map $\rho$ is constructed, we will show that it is 
an injection. Let $U_1$ and $U_2$ be two isotropic 
3-subspaces of $V$ such that $L_{U_1}=L_{U_2}$.
We can assume that $U_1={\rm Span}(v_1,v_2,v_3)$ and
representation $\hat V$ is written explicitly as above. It is 
clear that the only elements of $V$ that annihilate $L_{U_1}=\mathbb
C\hat  1$  are the ones from $U_1$, so $U_1=U_2$.

To calculate the image of $\rho$ observe that for any element $z$ 
of $\hat V$ its annihilator in $V$ is an isotropic subspace. Indeed,
an anticommutator 
of two elements of annihilator must be zero on $z$, but anticommutators
are given by the pairing. So the necessary and sufficient condition
on $z$ to come from an isotropic 3-space is to have the dimension
of its annihilator at least 3 (which will then be exactly 3 by
the above argument). We do this calculation explicitly. If
$$
z= \alpha_{\emptyset}\hat  1 + \alpha_{5} \hat v_5 +
\alpha_6 \hat v_6 + \alpha_7 \hat v_7 +\alpha_{56} \hat v_5\hat 
v_6 + \alpha_{67} \hat v_6\hat v_7
+\alpha_{57} \hat v_5\hat v_7 +\alpha_{567} \hat v_5\hat v_6\hat v_7 
$$
and $v=\sum_{i=1}^7 \beta_i v_i$ then $v\cdot z=0$ if and only if
$$
\begin{pmatrix}
-\alpha_7&\alpha_6&-\alpha_5&\frac 1{\sqrt 2}\alpha_\emptyset&0&0&0\\
\alpha_{57}&-\alpha_{56}&0&-\frac 1{\sqrt 2}\alpha_5&\alpha_\emptyset&0&0\\
\alpha_{67}&0&-\alpha_{56}&-\frac 1{\sqrt 2}\alpha_6&0&\alpha_\emptyset&0\\
0&\alpha_{67}&-\alpha_{57}&-\frac 1{\sqrt 2}\alpha_7&0&0&\alpha_\emptyset\\
-\alpha_{567}&0&0&\frac 1{\sqrt 2}\alpha_{56}&\alpha_6&-\alpha_5&0&\\
0&0&-\alpha_{567}&\frac 1{\sqrt 2}\alpha_{67}&0&\alpha_7&-\alpha_6\\
0&-\alpha_{567}&0&\frac 1{\sqrt 2}\alpha_{57}&\alpha_7&0&-\alpha_5\\
0&0&0&-\frac 1{\sqrt 2}\alpha_{567}&\alpha_{67}&-\alpha_{57}&\alpha_{56}\\
\end{pmatrix}
\begin{pmatrix}
\beta_1
\\\beta_2
\\\beta_3
\\\beta_4
\\\beta_5
\\\beta_6
\\\beta_7
\end{pmatrix}
=
\begin{pmatrix}
0\\0\\0\\0\\0\\0\\0\\0
\end{pmatrix}
.
$$
Hence $\mathbb P\mathbb Cz$ is in the image of $\rho$ if and only if the
rank of the above matrix is at most $4$. Two
$7\times 7$ determinants of this matrix obtained by removing the
first and the last row are
$$
D_1=-\frac 1{\sqrt 2}
\alpha_{567}(\alpha_\emptyset\alpha_{567}+\alpha_6\alpha_{57}-
\alpha_7\alpha_{56}-\alpha_{67}\alpha_5)^3,
$$
$$
D_8=\frac 1{\sqrt 2} 
\alpha_\emptyset(\alpha_\emptyset\alpha_{567}+\alpha_6\alpha_{57}-
\alpha_7\alpha_{56}-\alpha_{67}\alpha_5)^3.
$$
So the image of $\rho$ must be contained in the union of the hypersurface
\bean\label{explicitrho}
\alpha_\emptyset\alpha_{567}+\alpha_6\alpha_{57}-
\alpha_7\alpha_{56}-\alpha_{67}\alpha_5=0
\eean
and the 
subspace $\alpha_\emptyset=\alpha_{567}=0$.
By an easy dimension count, the dimension of ${\rm IG}(3,V)$ is $6$.
Hence, the image of $\rho$ is a hypersurface, so it must coincide 
with the smooth conic above. 
We also note that ${IG}(3,V)$ is smooth, since it is a variety 
with a transitive action of $SO(V)$, which implies that $\rho$ is an
embedding. Finally, this embedding is equivariant by construction and
\Ref{equiv}.
\end{proof}

\begin{rem}\label{conicistriv}
In view of Proposition \ref{symtwo} it is clear that the image of $\rho$ is
given by the dual of the one-dimensional subrepresentation
$\mathbb C$ of $Sym^2(\hat V)$.
Indeed, the image of $\rho$ is given by an element of $Sym^2(\hat V^*)$
which is invariant under the action of $Spin(V)$, perhaps up to a constant.
Therefore it is the dual of the trivial one-dimensional
subrepresentation of $Sym^2(\hat V)$. 
\end{rem}
\subsection{Spinor embedding and Pl\"ucker embedding}
Recall the spinor embedding
$\rho$ of the isotropic Grassmannian ${IG}(3,V)$ to $\Pb \hat V$
described in Theorem \ref{rhomap}. Let $\mc O(1)$ be the tautological
line bundle on  $\Pb \hat V$. We denote by $\mathcal O_{{IG}(3,V)}(1)$  the
line bundle on ${IG}(3,V)$ obtained  
as the pullback of the line bundle $\mc O(1)$ on $\mathbb P\hat V$ under
the embedding $\rho$. Let $\mathcal O_{{IG}(3,V)}(2)$ be the line bundle
obtained as a square of  $\mathcal O_{{IG}(3,V)}(1)$.
 
Another standard embedding of ${IG}(3,V)$ is the Pl\"ucker
embedding ${IG}(3,V)\to \mathbb P^{34}$ given by
$U\mapsto \C \Lambda^3(U)\in
\mathbb P\Lambda^3 V\cong \mathbb P^{34}$.

Recall that we have $V\cong V^*$ as $Spin(V)$-modules and
$Sym^2(V^*)=\C\oplus\La^3(V)$, see Proposition \ref{symtwo}.
The following proposition connects these two embeddings.
\begin{prop}\label{02}
The Pl\"ucker embedding is given by the 
full linear system of global sections of 
$\mathcal O_{{IG}(3,V)}(2)$  
.
\end{prop}

\begin{proof}
Let $\mc O$, $\mc O(-2)$ and $\mathcal O_{{IG}(3,V)}$  
be the sheaves of all functions on $\Pb \hat V$, 
functions on $\Pb \hat V$ vanishing on the image of $\rho$ and functions 
on ${IG}(3,V)$ respectively.
We have 
a short exact sequence of coherent sheaves on $\mathbb P \hat V$
$$
0\to \mathcal O(-2)\to \mathcal O\to \rho_*
\mathcal O_{{IG}(3,V)}\to 0.
$$
We tensor multiply it by $\mathcal O(2)$ and get
$$
0\to \mathcal O\to \mathcal O(2)\to \rho_*
\mathcal O_{{IG}(3,V)}(2)\to 0.
$$
We take the corresponding long exact sequence in cohomology.
Since $H^1(\mathbb P\hat V,\mathcal O)=0$
we see that the sections of $\mathcal O_{{IG}(3,V)}(2)$ are
given by the quotient of $Sym^2(\hat V^*)$ by the equation of 
the image of the spinor embedding. 
By Remark \ref{conicistriv} and 
Proposition \ref{symtwo}, this space is isomorphic as a $Spin(V)$
representation 
to $\Lambda^3 V\cong \Lambda^3 (V^*)$. 
As a result, the global sections of $\mathcal O_{{IG}(3,V)}(2)$ define 
a $Spin(V)$-equivariant map to $\mathbb P\Lambda^3 V$.

The image of ${IG}(3,V)$ is a six-dimensional orbit of $SO(V)$,
and the stabilizer of the point which corresponds to $U\subset V$
coincides with the stabilizer of $U$. One can see that 
the only fixed point of the stabilizer of $U$ in $\mathbb P\Lambda^3V$
is $\C \Lambda^3 U$, which means that the map
of the global sections of $\mathcal O(2)$ 
is precisely the Pl\"ucker embedding.
\end{proof}

\subsection{Properties of the image of the spinor embedding}
Choose a non-zero quadratic form $\hat Q$ on $\hat V$ which
vanishes on $\rho({IG}(3,V))$. We denote by 
$I\subset V$ and $\hat I\subset\hat V$, $\C\hat I=\rho({IG}(3,V))$ the
sets of isotropic vectors with respect to $Q$ and $\hat Q$ 
respectively. We denote by $\hat B$ the bilinear form 
associated to $\hat Q$. 
We will also make a choice of the action
of $V$ on $\hat V$ (as mentioned earlier these choices differ
by an overall sign only). We use calculations of 
Theorem \ref{rhomap} to prove some easy results about these
forms and the action.

\begin{prop}\label{skew}
For any $v\in V$, $p,q\in \hat V$ there holds
$\hat B(v\cdot p,q)=-\hat B(p,v\cdot q)$.
\end{prop}

\begin{proof}
By linearity, it is sufficient to show this statement for $v$ with
$Q(v)=-\frac 12$. Every such $v$ can be given by $v_4$ in some
basis as in Lemma \ref{spin model} and we 
use the corresponding model of $\hat V$.
It is then sufficient to check the statement for basic vectors
$p$ and $q$ in $\hat V$. Since the action of $v_4$
on the basis elements is diagonal with entries 
$(-1)^{\rm deg}\frac 1{\sqrt 2}$, we have 
$$\hat B(v\cdot p,q)=
(-1)^{{\rm deg} p-{\rm deg} q}\hat B(p,v\cdot q).$$ 
It remains to
observe from the explicit formula \eqref{explicitrho}
that $\hat B(p,q)=0$ unless $\deg(p)+\deg(q)=3$.
\end{proof}

\begin{cor}\label{vpp}
For all $v\in V$ and $p\in \hat V$ we have 
$\hat B(v\cdot p,p)=0$.
\end{cor}
 \begin{proof}
 $\hat B(v\cdot p,p)=-\hat B(p,v\cdot p)=-\hat B(v\cdot p,p).$
 \end{proof}

\begin{cor}\label{vpvp}
For all $v\in V$ and $p,q\in \hat V$ we have 
$\hat B(v\cdot p,v\cdot q)=Q(v)\hat B(p,q).$
\end{cor}

\begin{proof}
$\hat B(v\cdot p,v\cdot q)=-\hat B(v^2\cdot p,q)= -\hat
B(-Q(v)p,q)=Q(v)\hat B(p,q)$. 
\end{proof}

\begin{cor}\label{buvp} 
For all $p\in \hat V$ and $u,v\in V$ there 
holds
$\hat B(u\cdot p,v\cdot p)=B(u,v)\hat Q(p).$
\end{cor}

\begin{proof}
We have 
$\hat B(vu\cdot p,p)=-\hat B(u\cdot p,v\cdot p)=-\hat B(v\cdot p,u\cdot
p)=\hat B(uv\cdot p,p)$.
So $\hat B(u\cdot p,v\cdot p)=-\frac 12 \hat B((uv+vu)\cdot p,p)=
-\frac 12 \hat B(-B(u,v)p,p) = B(u,v)\hat  Q(p)$.
\end{proof}

\begin{prop}\label{lu}
Let $u\in I$ be  an isotropic vector acting on $\hat V$. 
Then we have $\on{Ker}(u) =
\on{Im} (u)\subset \hat I\subset \hat V$, $\on{dim}(\on{Ker}(u))=4$. Moreover,
$\Pb\on{Ker}(u)\in\Pb \hat V$
coincides with $\rho(\{U\in IG(3,V)\ |\ u\in U\})$.
\end{prop}
\begin{proof}
We complete $u$ to a basis of $V$ as in Lemma \ref{spin model} 
so that $u=v_1$. If we
pick the corresponding basis in $\hat V$ then the 
$u\cdot\hat V$ is the span of $\hat 1,\hat v_5,
\hat v_6,\hat v_5\hat v_6$ which is easily seen
to belong to $\hat I$. It is clearly the kernel of $u$.

Moreover if $U\in IG(3,V)$ and $\C p=\rho(U)$ then the space of elements
in $V$ which kill $p$ is exactly $U$. Therefore $u\cdot p=0$ if and
only if $u\in U$.
\end{proof}

\begin{prop}\label{orthcomp}
Let $p\in \hat V$ be a non-isotropic vector. Then 
$V\cdot p$ is the orthogonal complement of $p$ in $\hat V$ 
under $\hat Q$.
\end{prop}

\begin{proof}
By Corollary \ref{vpp} the space $V\cdot p$ is orthogonal to $p$.
So we need to show that $V\cdot p$ has dimension $7$. 
Suppose $v\cdot p=0$ for some $v$. Then $Q(v)p=-v^2\cdot p=0$ so
$v\in I$. By Proposition \ref{lu}, $p\in \on{Ker}(v)\subset \hat I$,
which contradicts the condition that $p$ is non-isotropic.
\end{proof}

\subsection{Invariant surjections}

A linear surjective map $\psi:\hat V\to V$ is called an {\it invariant
  surjection} if there exists $p\in \hat V$ such that
  $\psi(p)=0$ and $\psi(v\cdot p)=v$.

\begin{lem}
If $\psi$ be an invariant surjection then $p$ is non-isotropic.
\end{lem}
\begin{proof}
If $p$ was isotropic then we would have a non-zero $u\in V$ such
that $u\cdot p=0$ which gives a contradiction: $0=\psi(u\cdot p)=u$.  
\end{proof}

Let $V_1,V_2$ be linear spaces.
We call a map of projective spaces $\Pb\pi: \Pb V_1\to \Pb V_2$ {\it a
projection} if it is induced from a surjective linear map of linear
spaces $\pi:\ V_1\to V_2$.

Let $\psi:\ \hat V\to V$ be an invariant surjection and let $\Pb\psi:\
\Pb\hat V\to \Pb V$ be the induced projection. Let $\mu_\psi:\
IG(3,V)\to \Pb V$ be the composition of spinor embedding $\rho:\
IG(3,V)\to \Pb \hat V$ and $\Pb\psi$.

\begin{prop}\label{2 to 1}
For any isotropic $v\in V$ the equation $\mu_\psi(U)=\C v$ has a unique
solution $U\in IG(3,V)$. Moreover $v\in U$. 

For any non-isotropic $v\in V$ the equation $\mu_\psi(U)=\C v$ has exactly
two  solutions $U_1, U_2\in IG(3,V)$. Moreover $U_1+U_2=U_1\oplus
U_2=v^\perp$. 
\end{prop}
\begin{proof}
Recall that the image of the spinor embedding is the set of isotropic
vectors $\hat I\subset \hat V$.

Let $v$ be a non-zero element in $V$. The preimage of $\mathbb Cv$
in $\hat V$ is spanned by $p$ and $v\cdot p$. We have $\hat B(p,v\cdot
p)=0$ by Corollary \ref{vpp}.  
A linear combination $\al p+\beta v\cdot p$ is in $\hat I$ if and
only if (we use Corollary \ref{vpvp})
\bean\label{ox}
\hat Q(\alpha p +\beta v\cdot p) = \hat Q(p,p) (\alpha^2 + Q(v) \beta^2)=0.
\eean

If $v$ is isotropic we have $\al=0$ and therefore there is a unique
space $U\in IG(3,V)$ such that $\mu(U)=\C v$.
Moreover, since $v\cdot p$ is in $Ker(v)$, we get $v\in U$ by
Proposition \ref{lu}.   
 
If $v$ is not isotropic, then 
there are exactly two linear combinations (up to a scalar),
satisfying \Ref{ox} which give two spaces $U_1,U_2\in IG(3,V)$ such that
$\mu(U_1)=\mu(U_2)=\C v$. Observe that $U_1\cap U_2=0$. Indeed, if
$u\in U_1\cap U_2$ then $u$ kills both $\rho(U_1)$ and $\rho
(U_2)$ and therefore kills $p$ which is a contradiction.

Finally let $u\in U_1$. Then for $\al,\beta$ as above, we have $u\cdot(\al p
+\beta v\cdot p)=0$. It follows that $(\al-\beta v)u\cdot
p=B(u,v)p$. The left hand side is isotropic in the sense of
$\hat Q$, because $\hat Q(u\cdot p)=Q(u)\hat Q(p)=0$,
$\hat Q(vu\cdot p) = Q(u)Q(v)\hat Q(p)=0$ and 
$\hat B(u\cdot p, vu\cdot p)=0$ by Corollary \ref{vpp}.
Since $p$ is non-isotropic, we obtain $B(u,v)=0$. Therefore 
$v\perp U_1$. Similarly $v\perp U_2$.
\end{proof}

The following theorem provides a criterion for invariant surjections.
\begin{theorem}\label{g2proj}
Let $\psi:\hat V\to V$ be a surjective map with kernel
spanned by a non-isotropic $p\in\hat V$. Then $\psi$ is an invariant
surjection if and only if for a generic $u\in V$ there exists 
an isotropic dimension $3$ subspace $U\subset V$ such that $u\perp U$
and $\mu_\psi(U)=\C u$.
\end{theorem}

\begin{proof}
The only if part follows from Proposition \ref{2 to 1}.

Now we show the if part of the theorem. 
Consider the subset $Z$ of $IG(3,V)$ which consists 
of $U$ such that $\mu_\psi(U)\perp U$. Clearly, this is 
an algebraic subset. On the other hand, it has dimension
at least $6$, so the irreducibility of $IG(3,V)$ implies that 
$Z=IG(3,V)$ and 
$\mu_\psi(U)\perp U$ for all $U$. In particular,
if $u=\mu_\psi(U)$ is isotropic then $u\in U$,
since there are no isotropic $4$-subspaces of $V$.
This implies $u\cdot \rho(U)=0$.

The preimage of $u$ under $\psi$ contains $p$ and
$\rho(U)$. Since $u\cdot p\neq 0$, the only solution
up to scaling of 
$u\cdot \rho(q)=0$ for $\C\psi(q)=\C u$ is $\rho(U)$.
This implies that $\rho(U)\perp p$, since otherwise 
there would be two such solutions. This characterizes
$\rho(U)$ as the intersection of $\psi^{-1}\C u$ 
and $V\cdot p$.

Denote by $g:V\to V$ the map $g(v)=\psi(v\cdot p)$.
Surjectivity of $\psi$ and Proposition \ref{orthcomp} imply
that $g$ is invertible, and we denote its inverse by $h$.
We have $\psi(h(u)\cdot p) = u$, which implies that 
$\C h(u)\cdot p = \rho(U)$ and $uh(u)\cdot p =0$
for all isotropic $u\in I$.

Let $q_1,\dots, q_8$ be some basis in $\hat V$.
Write $uh(u)\cdot p
=\sum_{i=1}^8 a_i(u)q_i$ where $a_i(u)\in\C$.
For every linear function $r:\hat V\to \mathbb C$,
we have a quadratic function $r(uh(u)\cdot p)$ on $V$
which vanishes on $I$. Hence, it is proportional to $Q$. In particular
$a_i(u)=c_iQ(u)$ for some $c_i\in\C$. Therefore
there exists an element $p_1=\sum_{i=1}^8 c_iq_i\in \hat V$ such
that for all $u\in V$
$$
uh(u)\cdot p=Q(u)p_1=-u^2\cdot p_1.
$$
Since for generic $u$ its action on $\hat V$ has no kernel,
we get 
\bean\label{k}
h(u)\cdot p=-u\cdot p_1
\eean
for all $u\in V$.
This implies that $V\cdot p=V\cdot p_1$. By passing to 
the orthogonal complement, we get $p_1=cp$,
which means that $h(u)$ is a multiple of the identity.
We can scale $p$ to get $h(u)=u$, which implies that $\psi$ 
is an invariant surjection.
\end{proof}

We will need the following technical lemma later.
\begin{lem}\label{oneorbit}
For any $c\in \C^*$, the set
$\{ v\in \hat V\ |\ \hat Q(v)=c\}$ is an orbit of $Spin(V)$. 
\end{lem}

\begin{proof}
It is enough to show that any two points $p_1$, $p_2$ in
$\Pb\hat V$ lie in the same orbit of $Spin(V)$. 
Indeed, then every vector can be
translated to a multiple of another vector. If vectors have 
the same length, this multiple is $\pm 1$, and the lemma follows since 
$-{\bf id}\in Spin(V)$.

Let $p_1$ and $p_2$ be two points in $\Pb\hat V-\hat I$.
Draw generic lines $l_i$ through $p_i$ so that the line $p_i$ intersects
the conic $\hat I$ at two distinct points $q_{i,1}$ and $q_{i,2}$.
These points correspond to isotropic vector spaces $U_{i,1}$ and
$U_{i,2}$. 

Proposition \ref{2 to 1}, applied to an invariant surjection
 corresponding to point on $V$ which is on the line $p_i$,
assures that $U_{i,1}$ does not intersect $U_{i,2}$. 

Therefore there exists
an element of $SO(V)$ which maps $U_{1,1}$ to $U_{2,1}$ 
and $U_{1,2}$ to $U_{2,2}$. The corresponding elements of $Spin(V)$ 
will map $p_1$ to some point $p_3$ on the line $l_2$. It remains
to observe that the elements of $SO(V)$ that multiply all vectors
in $U_{2,1}$ by $\lambda$, vectors in $U_{2,2}$ by $\lambda^{-1}$
and fix the orthogonal complement $(U_1\oplus U_2)^\perp$ 
act transitively on the non-isotropic
points in $l_2$. Indeed, such elements scale the corresponding vectors in
$\hat V$ with coefficients $\pm\sqrt \lambda$ and $\pm\sqrt{ \lambda^{-1}}$
respectively.
\end{proof}

\subsection{A 3-form}\label{3-form sub}
Let $\psi:\ \hat V\to V$ be an invariant surjection.
In this section we define and study the trilinear skew-symmetric 3-form
that corresponds to $\psi$.

Let non-isotropic $p\in \hat V$ span the kernel of $\psi$.
We define a trilinear from $w_\psi$ by:
\be
w_\psi(a,b,c)=\hat B(abc\cdot p, p)(\hat Q(p))^{-1}
\ee
for all $a,b,c\in V$.
Note that the definition of $w_\psi$ does not depend on
the choice of non-zero $p$ in the kernel of $\psi$.

\begin{prop}\label{abcp}
The form $w_\psi$ is a trilinear skew-symmetric 3-form.
\end{prop} 
\begin{proof}
We check the skew-symmetry  (we use relations in Clifford algebra and 
Corollary \ref{vpp}):
\bea
\hat B(abc \cdot p, p) + \hat B(bac \cdot p, p) = - \hat B( B(a,b) c
\cdot p, p) = 0,\\
\hat B(abc \cdot p, p) + \hat B(acb \cdot p, p) = - \hat B( B(b,c) a
\cdot p, p) = 0.
\eea
\end{proof}

\begin{lem}\label{kernel dim 3}
Let $v\in V$ be an isotropic vector.  
Let $U$ be the kernel of the skew-symmetric 2-form
$w_\psi(v,\cdot,\cdot)$. Then $U$ is the unique 3-dimensional
isotropic space such that $\mu_\psi(U)=\C v$. In particular, $v\in U$.
\end{lem}
\begin{proof}
By Proposition \ref{2 to 1} the space $U\in IG(3,V)$ such
that $\rho(U)=\C v\cdot p$ is unique. Moreover $v\in U$. 

By the definition of $\rho$ and by the fact that the subspace of $V$
which kills the isotropic vector $v\cdot p\in\hat V$ is 3-dimensional, 
we obtain $u\in U$ if and only if $uv\cdot p=0$. 

If $u\in U$ then for any $a\in V$, we get
$$
w_\psi(a,u,v)=\hat B(auv\cdot p, p) = \hat B (0,p) = 0.
$$
Therefore $U\subset \on{Ker}\ w_\psi(v,\cdot,\cdot)$.

If for all $a\in V$, $w_\psi(a,u,v)=0$ then by Proposition \ref{skew}
\be
0=w_\psi (a,u,v)=\hat B(auv\cdot p,p)=-\hat B(uv\cdot p,a\cdot p).
\ee
But $V\cdot p=p^\perp$ by Proposition \ref{orthcomp}, hence
$uv\cdot p = \al p$ for some $\al\in \C$. But $\hat Q(uv\cdot p)=
Q(u)Q(v) \hat Q(p) = 0$, therefore $\al = 0$, and $uv\cdot p = 0$.
It follows that $\on{Ker}\ w_\psi(v,\cdot,\cdot)\subset U$.
\end{proof}

\begin{lem}\label{sympl}
Let $v\in V$ be a non-isotropic vector.  
Then the kernel of the skew-symmetric $2$-form 
$w_\psi(v,\cdot,\cdot)$ is one-dimensional and is 
spanned by $v$. 
\end{lem}

\begin{proof}
Suppose there is $u\in V$
such that for all $u_1\in V$ we have 
$$
\hat B(u_1vu\cdot p,p)=0.
$$
We rewrite this as $\hat B(vu\cdot p,u_1\cdot p)=0$,
which implies $vu\cdot p=\alpha p$ for some $\alpha\in\C$.
Acting on both sides by $v$, we obtain $-Q(v)u\cdot p = \alpha v\cdot
p$. Since $p$ is non-isotropic we have
$-Q(v)u=\alpha v$.
\end{proof}

\section{Group $G_2$}\label{G2 sec}
This section collects some facts about $G_2$. We believe
that  most of these facts
are standard and known. However in some cases we failed to find an adequate
reference and therefore we provide the proofs. 

\subsection{Definition of $G_2$}
Let $g_2$ be the semisimple complex Lie algebra with Cartan matrix
\be
\left(\begin{matrix}
~2 & -3~ \\
-1~ & ~2 
\end{matrix}\right).
\ee
We have $\dim g_2=14$. 

To describe the irreducible modules we choose a
Cartan subalgebra spanned by coroots $\al_1^\vee,\al_2^\vee$.
For $m,n\geq 0$, denote 
$L_{m,n}$ the irreducible highest module of highest weight
$m\om_1+n\om_2$\, where $\om_i$ are fundamental weights $\langle
\om_i,\al_j^\vee\rangle=\delta_i^j$. All irreducible modules have
such form. The dimension of $L_{m,n}$ is given by the Weyl formula:
\be
\dim L_{m,n}=\frac1{120}(m+1)(n+1)(m+n+2)(m+2n+3)(m+3n+4)(2m+3n+5).
\ee
In particular $\dim L_{0,0}=1$, $\dim L_{1,0}=7$, $\dim L_{0,1}=14$,
$\dim L_{2,0}=27$. We call the representations $L_{0,0}$ and $L_{1,0}$
trivial and vector representations respectively.

\medskip

Let $G_2$ be the connected, simply-connected complex Lie group with Lie
algebra $g_2$. 

Let $V$ be the 7-dimensional space with a non-degenerate 
bilinear symmetric form $B$ as before.
Let $\hat V$ be the spinor representation of $Spin(V)$.
The following fact is a $\C$-analog of  \cite[Theorem 5.5]{A}.

\begin{lem} Let $p\in \hat V$ be non-isotropic.
Then $G_2$ is isomorphic to the subgroup of $Spin(V)$ which fixes
$v$. 
\end{lem}

\begin{proof}
The stabilizer of $p$ is unaffected by scaling $p$.
Hence, 
by Lemma \ref{oneorbit}, the statement of this lemma is sufficient to
check for any non-isotropic $p$.

In view of \cite[Theorem 5.5]{A}, the Lie algebra of the compact
\emph{real} Lie group $G_2$ sits inside that of \emph{real} $Spin(7)$.
Consequently, the same is true for complex Lie algebras which 
by exponentiation implies that there is a map $G_2\to Spin(V)$
with a finite kernel. The determinant of the Cartan matrix of
$g_2$ is 1 and therefore the root and the weight lattices
of $g_2$ coincide. Then by \cite[Chapter 3, section 2.5]{GOV},
$G_2$ does not have finite subgroups and our map is injective.

Since $\hat V$ is a representation of $Spin(V)$,
it must split into representations of $G_2$, which can only 
be one copy of a trivial and one copy of a $7$-dimensional representation.
Hence there is an element $p\in \hat V$ whose stabilizer in $Spin(V)$
contains $G_2$. If $p$ was isotropic, then the Lie algebra $g_2$ would
act non-trivially on the tangent space to $\C p$ inside $\Pb \hat I$,
which is six-dimensional. This is impossible, which implies that 
$p$ is non-isotropic. 

It remains to show that the stabilizer $G$ of $p$ in $Spin(V)$
can not be bigger than $G_2$.
In view of the dimension count this can only happen
if $G$ is disconnected. Denote by $X$ the orbit of $p$, which is
isomorphic to the affine conic $\hat Q(*)=\hat Q(p)$ by 
Lemma \ref{oneorbit}. The group $Spin(V)$ is isomorphic to 
a fibration over $X$ with fiber $G$. The exact sequence of homotopy
groups gives an exact sequence $\pi_1(X)\to \pi_0(G)\to\pi_0(Spin(V))$,
so it is enough to check that $\pi_1(X)$ is trivial.

The variety $X$ is an orbit of $SO(\hat V)\cong SO(8)$. 
The stabilizer
of $p$ for the $SO(\hat V)$ action is isomorphic to $SO(\C p^\perp)
\cong SO(7)$. We have 
$$\pi_1(SO(7))\to \pi_1(SO(8))\to \pi_1(X)\to 1.
$$
The first map is an isomorphism, since it is induced by the 
embedding $Spin(7)\to Spin(8)$. This shows that $\pi_1(X)$ is trivial
and finishes the argument.
\end{proof}

We also use the projective version of this statement.
\begin{lem}\label{stabisG2}
The subgroup $G$ of $Spin(V)$ that fixes the point
$\mathbb Cp$ in $\Pb \hat V$ is disconnected.
The connected component of identity has index 2 and is 
isomorphic to the group $G_2$.
\end{lem}

\begin{proof}
Suppose that $g\in Spin(V)$ 
satisfies $g\mathbb Cp=\mathbb Cp$.
This implies $gp=cp$ for some constant $c$. Since $g$ preserves $\hat Q$, we
get $c=\pm 1$. The connected component of identity is characterized
by $c=1$ and therefore consists of elements of $Spin(V)$ that
fix $p$. On the other hand, $(-id_{\hat V})$ is an element
of $Spin(V)$ that has $c=-1$, which implies that the index of 
the $c=1$ subgroup is precisely two.
\end{proof}

We remark that in terms of the \emph{projective} action of
$SO(V)$ on $\mathbb P\hat V$ the isotropy subgroup of the point 
$\mathbb Cp$ is isomorphic to $G_2$. This is because
the map $Spin(V)\to SO(V)$ has kernel $\{\pm 1\}$. 
In particular $G_2$ is a subgroup of $SO(V)$ and acts on
$V$. Therefore $V$ is the vector representation of $G_2$. 

\subsection{$G_2$ and 3-forms}

We choose an identification $\Lambda^7 V\cong \C$.
A skew-symmetric 3-form $w$ on $V$ and a bilinear symmetric 2-form $b$
on the 7-dimensional space $V$ 
are called {\it associated} if there exists a non-zero constant
$c\in\C$, such that
for every two vectors $v_1,v_2\in V$ we have 
$b(v_1,v_2) = c w(v_1,*,*)\wedge w(v_2,*,*)\wedge w(*,*,*)$. 

A skew-symmetric 3-form on $V$ is called {\it non-degenerate} if it is
associated to a non-degenerate 2-form.

\begin{lem}\label{only inv}
Let $V$ be a vector representation of $G_2$. Then 

\begin{itemize}
\item
The space of $G_2$-invariant skew-symmetric
three-forms on $V$  is one-dimensional. Any such form is non-degenerate.
  
\item
The space of $G_2$-invariant symmetric bilinear 2-forms on $V$ is 
one-dimensional. Any such form is non-degenerate.

\end{itemize}

\noindent
In particular, any non-zero $G_2$-invariant skew-symmetric
three-form on $V$ is associated to any $G_2$-invariant 
symmetric bilinear 2-form.
\end{lem}
\begin{proof}
We have the following decompositions of $G_2$-modules:
\be
Sym^2(V)\cong L_{2,0}\oplus L_{0,0},\qquad \La^3(V)\cong L_{2,0}\oplus
L_{1,0}\oplus L_{0,0}.
\ee
These decompositions are easy to establish using the known
dimensions of $L_{m,n}$ and the weight decompositions of each
module.

It follows that the spaces of symmetric 2-forms and skew-symmetric
3-forms are one-dimensional. The bilinear form associated to a
$G_2$-invariant 3-form is $G_2$-invariant. 

So, we only have to show the $G_2$-invariant skew-symmetric form is non-degenerate. Since the bilinear form $b$ is a 
multiple of the non-degenerate form $B$,
 we simply need to show that the bilinear form associated to the $G_2$-invariant skew-symmetric $3$-form
is non-zero.  
We use the multiple of the $3$-form of Proposition \ref{abcp} given by 
$w(a,b,c)=\hat B(abc\cdot p, p)$ where $G_2=G_2(p)$ is realized as 
the subgroup of $Spin(V)$ that fixes $p\in \hat V$.

Consider a non-isotropic vector $v\in V$ and the
skew-symmetric $6$-form in $V/\C v$ given by
$w(v,*,*)\wedge w(v,*,*)\wedge w(v,*,*)$. It is sufficient
to show that this $6$-form is non-zero, 
since this will show that $b(v,v)\neq 0$.
By Lemma \ref{sympl} the form $w(v,*,*)$ is non-degenerate, so
its cube is non-zero.
\end{proof}

Let $\psi_p$ be an invariant surjection from a non-isotropic
vector $p$ and let  $G_2(p)$ be the group preserving $p$.
Recall that we constructed  a skew-symmetric trilinear form 
$w_\psi$ on $V$ in Section \ref{3-form sub}. 
\begin{prop}\label{threeform}
The form $w_\psi$ is a non-zero $G_2(p)$-invariant
3-form. In particular $w_\psi$ is associated to $B$.
\end{prop}

\begin{proof}

The form $w_\psi$ is $G_2(p)$-invariant because $p$ and $\hat B$ are
$G_2(p)$-invariant. The form $w_\psi$ is clearly non-zero in view of
Proposition \ref{orthcomp}.

The rest follows from Lemma \ref{only inv}.
\end{proof}

\begin{lem}\label{st=G2}
For a $G_2$-invariant 
non-degenerate skew-symmetric three-form $w$ the 
subgroup $G$ of $SL(V)$ that preserves $w$ is equal to
$G_2$.
\end{lem}

\begin{proof}
Clearly, $G\supseteq G_2$. Since $G$ preserves the associated
non-degenerate symmetric bilinear form $B$, it lies in $SO(V)$.

We can assume that $w$ is given by $\hat B(abc\cdot p,p)$.
By Lemma \ref{kernel dim 3} the action of $g$ on $\hat V$ 
preserves the set of $\on{Ker} w(v,*,*)$ for isotropic $v\in V$.
Hence it preserves the linear span $V\cdot p$
of this set.  It then preserves the orthogonal complement
of $V\cdot p$ which is the span of $p$ by Proposition \ref{orthcomp}. Lemma \ref{stabisG2} then finishes the proof.
\end{proof}

\begin{cor}
The $SL(V)$ orbit of a  $G_2$-invariant non-degenerate skew-symmetric three-form on $V$ 
is a dense open subset in the projective space of all non-zero skew-symmetric  three-forms on $V$ up to scaling.
\end{cor}

\begin{proof}
The dimension of 
the $SL(V)$ orbit of $w$ is equal to $\dim SL(V)-\dim G_2=
48-14=34=\dim \Pb \Lambda^3(V^*)$, which shows that the 
orbit is dense.
\end{proof}

For the sake of the completeness of exposition we prove a stronger
fact, which will not be used in the rest of the paper.
Namely, we now show that every non-degenerate skew-symmetric
three-form on $V$ lies in this  $SL(V)$ orbit, up  
to scaling. 
This implies that every non-degenerate skew-symmetric three-form
 is $G_2$-invariant for some $G_2\subset SL(V)$.
\begin{prop}\label{stab=G2}
Let $w$ be a trilinear skew-symmetric 3-form on $V$ associated to
$B$. Then the subgroup of $SL(V)$ which preserves $w$ is isomorphic to
$G_2$. 
\end{prop}
\begin{proof}
Let $w$ be a non-degenerate skew-symmetric three-form
and let $B$ be the associated non-degenerate bilinear
symmetric form. Clearly, the group $G$ that preserves $w$
must preserve $B$ as well, so it must lie in $SO(V)$. 

The argument of the last paragraph of the proof of 
Lemma \ref{only inv} shows that 
for every element $v\in V$ the bilinear form $B(v,v)$
is zero if and only if skew-symmetric bilinear form
$w(v,*,*)$ on $V/\C v$ has zero kernel. This implies 
that if $w(u,v,*)$ is zero then $u$ and $v$ are both
isotropic. Therefore, for any $u\in I$, the kernel 
of $w(u,*,*)$ consists of isotropic vectors of $B$ and
hence has dimension at most three. Consequently,
the dimension of the kernel of $w(u,*,*)$ is 
three if $u$ is isotropic and one if $u$ is not. In both
cases the kernel contains $u$.

The above discussion produces a map 
$\phi: \Pb I\to \Pb \hat I$ by mapping the isotropic vector
to the corresponding point in $IG(3,V)$ which is identified
with $\Pb\hat I$ via the spinor embedding. 

We will now show
that $\phi$ is an embedding. Suppose two 
linearly independent vectors $u_1,u_2$ have
the same kernel $U$ of $w(u_i,*,*)$. Then for a generic
element $u\in U$ the kernel of $w(u, *,*)$ contains 
$u$, $u_1$ and $u_2$, which implies that $w(u,u',*)=0$
for all $u,u'\in U$. Pick two vectors $v\in V$ and $u\in U$
with $B(u,v)\neq 0$. If we look
at 
$$
B(u,v)=w(u,*,*)\wedge w(v,*,*)\wedge w(*,*,*)
$$
we see that the right hand side is zero. Indeed, if we were
to evaluate on the basis of $V$ that extends the basis
$u_1,u_2,u_3$ of $U$, then in each term either we have
$w(u,u_i,*)$ or $w(v,u_i,u_j)$, or $w(u_i,u_j,*)$. Contradiction.
We can similarly show that the map $\mu$ does not kill
tangent vectors. Indeed, if we 
have $\on{Ker} (u_1+\epsilon u_2)=U\mod \epsilon^2$,
we have $w(u_1,U,*)=w(u_2,U,*)=0$, and the previous argument
works.

As a result, the image $\phi(\Pb I)$ of $\phi$ is a smooth conic of
dimension $5$ inside a smooth conic of dimension $6$ in 
a projective space of dimension $7$. By strong Lefschetz
theorem, $\phi(\Pb I)$ is a complete intersection of $\Pb \hat I$
and a hypersurface in $\Pb \hat V$. It is easy to see 
that this hypersurface must in fact be a hyperplane, by 
comparing numerical invariants of complete intersections
and conics. Hence, $\phi(\Pb I)$ spans a hyperplane in
$\Pb\hat V$. Denote by $V_1$ the corresponding codimension
one subspace in $\hat V$. The restriction
 of $\hat B$ to
$V_1$ gives a smooth conic $\phi(\Pb I)$, which implies that
the orthogonal complement to $V_1$ is non-isotropic. 

Every element $g\in G$ must preserve
$\phi(\Pb I)$ as a set. Hence it preserves $V_1$ and 
its orthogonal complement, so it preserves
a point in $\Pb \hat V$. Lemma \ref{stabisG2} then 
shows that $G$ is contained in $G_2$ so it has dimension
at most $14$. Then the $SL(V)$ orbit of $w$ is 
dense in $\Pb \Lambda^3(V^*)$ so it must coincide with
the orbit of a $G_2$-invariant form. This shows that
$G$ is a conjugate of some $G_2\subset SL(V)$.
\end{proof}

Note that $SO(V)$ acts transitively on 3-forms associated to $B$
(considered up to a constant). 
Indeed, $GL(V)$ acts transitively on the set of all non-degenerate
3-forms (because each orbit is an algebraic variety of dimension 35). And any element
which maps a 3-form associated to $B$ to a 3-form associated to $B$
has to preserve $B$ (all up to a constant). 

In particular the dimension of the variety of 3-forms (considered up
to a constant) associated to $B$ is 7.

It follows that all subgroups of $SO(V)$ isomorphic to $G_2$ are
conjugated. We have a bijection between the set of subgroups of $SO(V)$
isomorphic to $G_2$ and 3-forms associated to $B$ (considered up to a
constant). We also have a bijection between the set of subgroups of $SO(V)$
isomorphic to $G_2$ and $\Pb \hat V-\hat I$.



\subsection{The homogeneous space $G/B$}
Recall that Borel subgroup of a semisimple Lie group is a maximal
connected solvable subgroup. 
We give explicit description of flag varieties $G/B$ of types $G_2,C_3,A_6$.

A full flag $F$ in $V$ is the chain 
$
F_1\subset  F_2 \subset  F_3 \subset  F_4 \subset  F_5
\subset  F_6,
$
where $ F_i\subset V$ are subspaces of dimension $i$. 

\medskip

Let $\mc F$ be the set of all full flags in $V$. The group $SL(V)$ acts
transitively on $\mc F$, this action identifies $\mc F$ with the
quotient of $SL(V)$ by a Borel subgroup. Then $\mc F$ is an
algebraic variety of dimension $21$ which is called the flag variety
of $A_6$ type. 

\medskip

A full flag $F$ is called isotropic if $B(F_i,F_{7-i})=0$.
Let $\mc F^\perp\subset \mc F$ be the set of all isotropic flags in
$V$. The group $SO(V)$ acts transitively on $\mc F^\perp$, this action 
identifies $\mc F^\perp$ with the quotient of $SO(V)$ by a Borel subgroup. 
Then $\mc F^\perp$ is an algebraic variety of dimension $9$ 
which is called the flag variety of $B_3$ type. 

\medskip
We choose a 3-form $w$ on $V$ associated to $B$ and $G_2(w)$ the
subgroup of $SO(V)$ which preserves $w$.
An isotropic flag $F$ is called $G_2$-isotropic 
if $F_3=\on{Ker}w(F_1,*,* )$.
Let $\mc F^{\prp}\subset \mc F^\perp $ be the set of all
$G_2$-isotropic flags in $V$. 

\begin{prop}\label{G/B}
The set $\mc F^{\prp}$ is a smooth algebraic variety of dimension $6$..
The group $G_2(w)$ acts transitively on the space $\mc F^{\prp}$, this action 
identifies $\mc F^{\prp}$ with the quotient of  
$G_2(w)$ by a Borel subgroup. 
\end{prop}
\begin{proof}
 Elements 
of $\mc F^{\prp}$ can be identified with 
flags $F_1\subset F_2$ such that for bases $u$
and $(u,v)$ of $F_1$ and $F_2$ respectively
there holds $B(u,u)=B(u,v)=B(v,v)=0$ and $w(u,v,*)=0$. 
Consequently, $\mc F^{\prp}$ is a smooth
variety of dimension $6$ which admits a fibration
$\mc F^{\prp}\to \Pb I$ to a smooth conic of dimension $5$.
The fibers are isomorphic to $\Pb^1$. In fact, the fibration
is the projectivization of the rank two vector bundle whose fiber 
over $F_1$ is given by $\on{Ker} w(F_1,*,*)/F_1$.

Denote by $B_1$ the
stabilizer of $F_1\subset F_2$ in $\mc F^{\prp}$.
It is easy to see that it stabilizes a complete flag in $V$
and so lies in the Borel subgroup of $SL(V)$. Consequently,
$B_1$ is solvable which implies that its connected component
of identity is contained in a Borel subgroup $B_2$ of $G_2(w)$. 

The dimension of $B_2$
is $8$, which implies that the dimension of $B_1$ is at most
$8$. Since the dimension of $G_2$ is $14$ the dimension
of $G_2(w)$ orbit of the flag is at least $6$, hence it is exactly $6$.
Since $\mc F^{\prp}$  is irreducible, this means that
all $G_2(w)$ orbits on $\mc F^{\prp}$ are Zariski
dense hence
there is only one orbit. We have thus shown the transitivity
of the action.

The above argument also shows that the connected 
component of identity $B_1^\circ\subseteq B_1$ is equal to $B_2$,
since their Lie algebras have the same dimension and one 
is contained in the other.

We now claim that $B_1$ is in fact connected. If 
$B_1$ was not connected, the variety $\mc F^{\prp}$ 
would admit an unramified covering from the variety
$G_2(w)/B_1^\circ$. However, the fibration description of 
$\mc F^{\prp}$ shows that it is unirational,
i.e. admits a dominant morphism from a projective space.
Hence it has trivial fundamental 
group by \cite{S}. 
\end{proof}
The quotient of $G_2$ by a Borel subgroup is called 
the flag variety of $G_2$ type. Proposition \ref{G/B} says that $\mc
F^{\prp}$ is isomorphic to the flag variety of $G_2$ type.

\section{Self-self-dual spaces of polynomials}\label{ssd spaces}
In this section we 
define the self-self-dual spaces of polynomials and show that such
spaces possess a natural skew-symmetric 3-form.

\subsection{Self-dual spaces of polynomials of dimension 7}
Here we recall main facts about self-dual spaces of polynomials proved
in \cite{MV}.
  
Let $V\in\C[x]$ be a vector space of polynomials of dimension $7$. The
space $V$ has {\it a base point z} if all polynomials in $V$ vanish at $z$.
We always assume that $V$ has no base points.

Let $W(u_1,\dots,u_k)$ denote the Wronskian of functions $u_1,\dots,u_k$, 
\be 
W(u_1,\dots,u_k)=\on{det}(u_{i}^{(j-1)})_{i,j=1}^k.
\ee
Let $U_i$ be the (monic) greatest common divisor of the set of all
Wronskians $W(u_1,\dots,u_i)$
where $u_1,\dots,u_i\in V$. The following lemma is easy, cf. for
example Lemma 4.9  in \cite{MV1}.

\begin{lem}
There exist unique polynomials $T_1,\dots,T_6$ such that
$U_i=T_1^{i-1}T_2^{i-2}\dots T_{i-1}$. $\Box$
\end{lem}
We call polynomials $T_i$ {\it the ramification polynomials} of $V$.

For $u_1,\dots,u_i\in V$, we call {\it the divided Wronskian} the
polynomial defined by 
\be
W^\dagger(u_1,\dots,u_i)=\frac{W(u_1,\dots,u_i)}{U_i}=\frac{W(u_1,\dots,u_i)}{T_1^{i-1}T_2^{i-2}\dots T_{i-1}}.
\ee
Note that $W^\dagger(u_1,\dots,u_7)\in\C$.

Space $V$ is called {\it self-dual} if 
\be
V=\{W^\dagger(u_1,\dots,u_6)\ |\ u_1,\dots,u_6\in V \}.
\ee
If $V$ is self-dual then $T_i=T_{7-i}$.

By Corollary 6.5 in \cite{MV} a self-dual space of polynomials $V$
possesses a bilinear symmetric non-degenerate form given by
\bean\label{2form}
B(u_1,u_2)=W^\dagger(u_1,v_1,\dots,v_6),\qquad \on{if} \
u_2=W^\dagger(v_1,\dots, v_6).
\eean 

\begin{lem}\label{always a square} 
Let $V$ be a self-dual space of polynomials of dimension 7. 
Let $u_1,u_2,u_3\in V$ be such that $B(u_i,u_j)=0$. Then the divided
Wronskian $W^\dagger(u_1,u_2,u_3)$ is a square of a polynomial.
\end{lem}
\begin{proof}
By Lemma 6.16 in \cite{MV}, the divided Wronskian is a square for
generic isotropic 3-space. The condition of Lemma 6.16 in \cite{MV} is
satisfied for a generic isotropic 3-space by Theorem 7.10 in \cite{MV}. 
Then the divided Wronskian of any isotropic 3-space is a square.
\end{proof}

\subsection{Definition of self-self-dual spaces and first properties}
Let $V$ be a self-dual space of polynomials of dimension $7$.

\begin{defn}\label{selfself}
A self-dual space of polynomials $V$ of dimension 7 is called
self-self-dual if for a generic element
$f$ of $V$ there exists a three-dimensional isotropic subspace
$U\subset V$ such that $f\perp U$ and the divided
Wronskian of $U$ is (up to a constant) the
square of $f$. 
\end{defn}

\begin{rem} The space of polynomials of degree at most $6$ is self-self-dual, see Section \ref{monex}.
\end{rem}

Let $V$ be a self-self-dual space of polynomials. Then $V$ is
self-dual and all divided Wronskians of isotropic 3-spaces are squares
of polynomials by Lemma \ref{always a square}.
Denote by $\bar V$ the span of all square roots of all divided Wronskians
of three-dimensional isotropic subspaces $U\subset V$. Then we have a map
\be
\mu:\ { IG}(3,V)\to \mathbb P\bar V. 
\ee

\begin{prop}\label{mainproj}
There exists a projection $P:\ \Pb\hat V\to \Pb \bar V$ such that
the map $\mu$ is given by a composition of the spinor embedding
$\rho:\ {IG}(3,V)\to \Pb\hat V$ and projection $P$.
\end{prop}

\begin{proof}
Let $L$ be the tautological line bundle on $\bar V$.
Every map from an algebraic variety $M$ to a projective space $\bar V$ 
is given by sections of the line bundle on $M$ which is 
the pull-back of $L$. So, we describe the pullback bundle $\mu^*(L)$. 
To do that we consider the pullback of $L^2$.

Recall that there is a bijective correspondence between line bundles
and $\C^*$-bundles on any algebraic variety. We denote the
$\C^*$-bundle corresponding to a line bundle $M$ by $\tilde M$.

The square of $\tilde L^{-1}$ can be thought of as the variety 
of non-zero squares $f^2, f\in \bar V$ with the obvious map to
$\mathbb P \bar V$. 

Let $\tilde R$ 
be the subvariety of $\La^3 V$ which corresponds to isotropic 3-spaces with a choice of a volume form:
\be
\tilde R=\{u_1\wedge u_2\wedge u_3\ | u_1\wedge u_2\wedge u_3\neq 0, 
\ u_i\in V, \ B(u_i,u_j)=0\}\subset \La^3 V.
\ee 
Then $\tilde R$ is $\C^*$-bundle over ${IG}(3,V)$
which is the dual of the $\C^*$-bundle associated to the 
Pl\"ucker bundle.
By Proposition \ref{02}, $R$ is $\mathcal O_{{IG}(3,V)}(-2)$.

The map $\mu$ induces a map of $\C^*$-bundles 
from $\tilde R$ to the pullback of $\tilde L^{-2}$ 
which maps $u_1\wedge u_2\wedge u_3$ to $f^2$ if
$W^\dagger(u_1,u_2,u_3)=f^2$. Every map of $\C^*$ bundles which is
isomorphism on the bases is an isomorphism.
Therefore $\tilde R$ is the pullback of $\tilde L^{-2}$ and $\mathcal
O_{{IG}(3,V)}(2)$ is the pullback of $L^2$.

The line bundles of any  smooth hypersurface in a projective
space of dimension at least $4$
are integer powers of the pullback of the tautological bundle of the
projective space. This follows for example from the 
strong Lefschetz theorem \cite{GH}. In particular, all bundles
on ${IG}(3,V)$ are powers of $\mathcal O_{{IG}(3,V)}(1)$.

Consequently, the pullback of $L$ is the spinor line bundle
$\mathcal O_{{IG}(3,V)}(1)$. Thus the map $\mu$ is given by a subspace of
the global sections of $\mathcal O_{{IG}(3,V)}(1)$. The global
sections of this bundle are a subspace in $\hat V$. The lemma follows.
\end{proof}

\begin{prop}\label{there is a p}
The space $\bar V$ is equal to $V$. Moreover, there exists an
invariant surjection $\psi:\ \hat V\to V$, such that projection $P$ is
the projectivization of $\psi$.
\end{prop} 
\begin{proof}
The condition of Definition \ref{selfself}
implies that $\bar V$ contains $V$. 
Assume that it is bigger than $V$. By Proposition \ref{mainproj}
the dimension of $\bar V$ is at most the dimension of $\hat V$,
so it must be exactly $8$ and the projection of Proposition
\ref{mainproj} must in fact be an isomorphism. Then 
the image of $\mu$ is a non-singular conic in $\mathbb P\bar V$
and does not contain a generic element of $\mathbb P\bar V$.

Now we know that the map $\mu$ is a composition of the
embedding ${IG}(3,V)$ as a conic in $\mathbb P\hat V$ and some
projection given by a surjective map $\psi:\hat V\to V$.
The fact that the map $\mu$ is well-defined implies that the kernel
of $\psi$ does not lie in ${IG}(3,V)$. 
The condition of Definition \ref{selfself} implies that
$\psi$ satisfies the assumption of Theorem \ref{g2proj} and therefore $\psi $ is an invariant surjection.
\end{proof}

We have an immediate corollary. 
\begin{cor}\label{2 to 1 wr}
Every vector $f$ in $V$ is  
a constant times a square root of a divided Wronskian 
of a three-dimensional isotropic subspace $U$ of $V$. If $f$ is
isotropic, such $U$ is unique and $f\in U$. If $f$ is not isotropic,
there are exactly 
two such spaces $U_1$ and $U_2$, and $U_1+U_2=U_1\oplus U_2=f^\perp$.
\end{cor}
\begin{proof}
Follows from Proposition \ref{2 to 1}.
\end{proof}

\begin{lem}
If $V$ is a self-self-dual space then $T_1=T_3=T_4=T_6$ and $T_2=T_5$.
\end{lem}
\begin{proof}
We only need to show $T_1=T_3$. All $T_i$ have zeroes only at
$z_1,\dots,z_n$. Let us show that $T_1$ and $T_3$ have the same 
order of zero at $z_1$. 

Let $T_1$, $T_2$, $T_3$ have orders $t_1,t_2,t_3$ of zero 
at $z_1$ respectively.
Let $u_1,u_2,\dots,u_7$ be a basis of $V$ such that order of zero at
$z_1$ of $u_i$ is strictly smaller than the order of zero at $z_1$ of
$u_{i+1}$.  Then the orders of $u_i$ at $z_1$ are given by 
\be
0,t_1+1,t_1+t_2+2,t_1+t_2+t_3+3,t_1+t_2+2t_3+4,t_1+2t_2+2t_3+5,2t_1+2t_2+2t_3+6.
\ee  
Then $W^\dagger(u_5,u_6,u_7)$ has the maximal possible
order of zero at $z_1$. Since $V$ is self-self-dual, this implies 
$W(u_5,u_6,u_7)=cu_7^2T_1^2T_2$ for some $c\in\C^*$. 
But the order of zero of
$W(u_5,u_6,u_7)$ is $4t_1+5t_2+6t_3+12$ and the order of zero of 
$u_7^2T_1^2T_2$ is $6t_1+5t_2+4t_3+12$. It follows that $t_1=t_3$.

\end{proof}

\subsection{The 3-form}
Let $p\in \hat V$ span the kernel of the map  $\bar P$ which
corresponds to the projection $P$ in Proposition \ref{mainproj}. We
denote by $G_2(p)$ 
the subgroup of $SO(V)$ which preserves $p$. The group
$G_2(p)$ is isomorphic to $G_2$.
We characterize $G_2(p)$ as the group of linear 
transformations of $V$ that are compatible with the Wronskian structure.
\begin{prop}\label{g2v}
The group $G_2(p)$ is the set of all elements $g\in SO(V)$ such that  
for any basis $(a_1,a_2,a_3)$ of any isotropic 3-subspace 
of $V$, any $f\in V$ satisfying 
$W^\dagger(a_1,a_2,a_3) = f^2$ there holds
\begin{equation}\label{star}
W^\dagger(ga_1,ga_2,ga_3) = (gf)^2.
\end{equation}
\end{prop}

\begin{proof}
For each $g\in SO(V) $ 
consider a map $\mu_g:{IG}(3,V)\to \mathbb PV$ 
given as follows. For every isotropic 3-space $U$ pick a basis $(a_1,
a_2,a_3)$ and define
$$
\mu_g(U):= g^{-1}\Big(\sqrt {W^\dagger(g(a_1),g(a_2),g(a_3))}\Big).
$$
While the square root is defined up to a sign only, the result makes
sense as a point in $\mathbb P V$. 

The proof of Proposition \ref{mainproj} is applicable to the map
$\mu_g$ as well as to $\mu$. Consequently, the map is given by a
projection from a point $P_{\mu_g}$ in $\mathbb P\hat V$ which is uniquely 
determined by $g$. Moreover, the 
argument of Proposition \ref{mainproj} shows that the 
map is uniquely determined by 
that point, since Theorem \ref{g2proj} still applies.
Our construction is $SO(V)$-equivariant,
so we have $P_{\mu_g}=g\mathbb C p$ 
for the natural projective action 
of $SO(V)$ on $\mathbb P\hat V$. 
As a result, equation \eqref{star} implies that $g\in G_2(p)$.

On the other hand, for every $g\in G_2(p)$
we have $\mu_g=\mu$ which translates into 
$$
W^\dagger(g(a_1),g(a_2),g(a_3)) = c(a_1,a_2,a_3,g)(g(f))^2.
$$
for some constants $c\in \mathbb C^*$ for all choices of $a_i$
and $f$ as in the statement of the proposition. The constant 
$c$ depends only on $g$ and  the point $V_1\in 
{IG}(3,V)$. Since it is homomorphic and ${IG}(3,V)$ is 
compact, it must depend on $g$ only. 
Since $G_2$
has no non-trivial one-dimensional representations, 
we see that $c=1$, so elements of $G_2(p)$
satisfy \eqref{star}.
\end{proof}

Define a linear map $\nu^*:\ \La^3 V\to Sym^2 V^*$ by the formula
\be
\langle\nu^*(a\wedge b \wedge c),v\otimes v \rangle=
\hat B (abcv\cdot p,v\cdot p),
\ee
where $a,b,c,v,\in V$.
The argument of Proposition \ref{abcp} shows that $\nu^*$ is well defined.

Let $\nu:\La^3 V\to Sym^2 V$ be the linear map  which is a composition
of $\nu^*$ with the natural identification
of $Sym^2 V^*$ and $Sym^2 V$ induced by $B$.

We remark that when one scales $p$ or $\hat B$ one scales $\nu$.

\begin{prop}\label{phinu}
There exists a unique linear map 
$\phi:\La^3 V\to Sym^2 V$ such that
$\phi(a\wedge b\wedge c)  = f\otimes f$
for all isotropic pairwise orthogonal $a,b,c$ with $W^\dagger (a,b,c) =f^2$.
Moreover, there is a constant $C$ such that $\phi=C\nu$.
\end{prop}
\begin{proof}
We need to show that there is a constant $C$ such that
$\nu(v_1\wedge v_2\wedge v_3) =Cf\otimes f$ for all bases
$(v_1,v_2,v_3)$ of an isotropic subspace $U$ and $f$ such that 
$W^\dagger (v_1,v_2,v_3) = f^2$. 

We extend $(v_1,v_2,v_3)$ to a basis of $V$ as in Theorem \ref{rhomap}.
To calculate 
$\hat B(v_1v_2v_3 v \cdot p, v \cdot p)$, observe 
that $\hat B(v_1v_2v_3 \cdot q,q)$ is proportional 
to $\alpha_{567}^2$ in the notations of the proof of Theorem
\ref{rhomap}. This is also proportional to 
$(\hat B (\hat 1,q))^2$
Since we have $\hat 1 = \gamma_1 p + \gamma_2f\cdot p$
for some $\gamma_1$ and $\gamma_2$, we see that 
$$\langle\nu^*(a\wedge b\wedge c),v\otimes v\rangle =
 C(a,b,c) (\hat B(\gamma_1p,v\cdot p) + \hat B(\gamma_2f\cdot p ,
 v\cdot p))^2
$$
$$
=C_1(a,b,c) B(f,v)^2.
$$
Hence $\nu(a\wedge b\wedge c) = C_1(a,b,c) f\otimes f$
for some nonzero $C_1(a,b,c)$. 
We now observe that $C_1$ in fact depends on the point of 
$IG(3,V)$ only. Since it is clearly holomorphic, it must be
a constant. 

We now show the uniqueness of $\phi$. 
If we have $\phi_1$ and $\phi_2$ that satisfy the conditions
of the proposition their difference $\phi_3$ satisfies
$\phi_3(a\wedge b\wedge c) = 0$ for all bases $(a,b,c)$
of an isotropic subspace. Since the $SO(V)$ representation
$\La^3 V$ is irreducible, such $a\wedge b\wedge c$ span
all of $\La^3 V$, which implies that $\phi_3=0$.
\end{proof}

Let $m:\ Sym^2 V\to \C[x]$ be the multiplication map, sending
$f\otimes g + g\otimes f \mapsto 2fg$.
\begin{cor}\label{general wrons}
For all $a,b,c\in V$, we have $W^\dagger(a,b,c)=m(\phi(a\wedge b
\wedge c))$.
\end{cor}
\begin{proof}
By Proposition \ref{phinu}, the corollary is 
true if $a,b,c$ span an isotropic subspace. Such triples span
$\La^3V$, since they span an $SO(V)$-subrepresentation in $\La^3V$,
but $\La^3V$ is an irreducible $SO(V)$-module.
Then the generic case of the corollary follows from linearity. 

\end{proof}

Now we obtain the $G_2(p)$-invariant 3-form.
\begin{theorem}\label{wrform}
There is a unique skew-symmetric 3-form $w\in \Lambda^3V^*$ such that  
\bean\label{3fm}
w(a_1,a_2,a_3)=B(f,f)
\eean
for all pairwise orthogonal isotropic $a_1,a_2,a_3$ with
$W^\dagger(a_1,a_2,a_3)=f^2$.
An element $g\in SO(V)$ preserves the form $w$ if and only if $g\in
G_2(p)$. Finally, the form $w$ is associated to $B$.
\end{theorem}
\begin{proof}
We set $w(a,b,c)=B(\phi(a\wedge b\wedge c))$, where $\phi$ is the map
in Proposition \ref{phinu}. This form clearly satisfies all the
properties described in the theorem. In particular it is associated to
$B$ because all $G_2$-invariant forms are associated to $B$ by Lemma
\ref{only inv}.

\end{proof}

\begin{cor}
Let $v\in V$, $Q(v)=0$. Let $U$ be the kernel of the skew-symmetric 2-form
$w(v,\cdot,\cdot)$. Then $U$ is the unique 3-dimensional isotropic
3-space, $U\in IG(3,V)$ such that
$W^\dagger(U)$ equals (up to a constant) $v^2$. Moreover $v\in U$.  
\end{cor}
\begin{proof}
Follows from Lemma \ref{kernel dim 3}. 
\end{proof}

\section{Populations of critical points}\label{populations}
Populations of critical points are defined for any Kac-Moody
algebra in \cite{MV}. The motivation for this definition is algebraic
Bethe Ansatz method.
In this paper we study the populations of critical points associated with
$G_2$ root systems. We use properties of populations of
$A_6$ and $C_3$ type. We recall known facts about those cases first.


\subsection{$A_6$-populations}
Consider the root system of type $A_6$. 
We have 
\be
(\al_i,\al_i)=2, \qquad (\al_i,\al_{i\pm 1})=-1, 
\ee
and other scalar products are zero. The root system of type $A_6$
corresponds to Lie algebra $sl_7$. 
The Weyl group of $sl_7$ is generated by simple reflections
$s_i$. The shifted action of Weyl group on $A_6$-weights is
\be
s_i\cdot \la=\la-(\la+\rho,\al_i)\al_i,
\ee 
where $(\rho,\al_i)=1$ for all $i$.

Fix polynomials
$T_i(x)=\prod_{s=1}^n(x-z_s)^{m_s^{(i)}}$.

where $i=1,\dots,6$, 
$x$ is a formal variable and $z_i,m_s^{(i)}$ are parameters.

A $6$-tuple of polynomials $\bs y=(y_1,\dots,y_6)$ is called {\it
  generic} if $y_i$ have
   no multiple roots and if $y_i,y_{i+1}$ have no
  common roots. A generic $6$-tuple of polynomials $\bs y$ is called
  {\it a Bethe $6$-tuple} of  type $A_6$ 
 if there exist polynomials $\tilde
  y_1,\dots,\tilde y_6$ such that 
\be
W(y_i,\tilde y_i)=T_iy_{i-1}y_{i+1}.
\ee
Note that $\tilde y_i$ is unique up to addition of a scalar multiple of $y_i$.

In what follows we always consider $N$-tuples of polynomials up to
multiplication of its coordinates by non-zero scalars.

If $\bs y$ is a Bethe $6$-tuple and the $6$-tuple $\bs
y^{(i)}=(y_1,\dots,\tilde y_i, \dots,y_6)$ is generic then $y^{(i)}$
is also a Bethe $6$-tuple of $A_6$ type, see Theorem 3.7 in \cite{MV}.
We call $\bs y^{(i)}$ {\it an immediate descendent of $\bs y$ in the
direction $i$}.

Let $\bar P$ be the minimal set of Bethe $6$-tuples which contains
$\bs y$ and such that for all $\tilde{ \bs y}\in\bar P$ all immediate
descendants  $\tilde {\bs y}^{(i)}$ are also in $\bar P$. 

The degrees of all coordinates of all tuples of
$\bar P$ are simultaneously bounded.
{\it The $A_6$-population of Bethe $6$-tuples originated at $\bs y$} 
is the closure of $\bar P$ in $(\Pb (\C[x]))^6$.

\bigskip

The set of populations of Bethe $6$-tuples associated with polynomials
$T_i$ is in one-to-one correspondence with the set of $7$-dimensional
spaces of polynomials with no base points and ramification polynomials $T_i$.
Given a population $P$ and $\bs y\in P$, we define
\bean\label{Diff}
D & =&  (   \partial  - 
\ln'  ( \frac { \prod_{s=1}^6 T_s } { y_6 }  ) )
( \partial - \ln' ( \frac {y_6  \prod_{s=1}^{5} T_s} {y_{5} }  ) )
\dots ( \partial  -  \ln'  ( \frac { y_2 T_1 }{ y_1 }  ) ) 
( \partial - \ln' (  y_1 ) ) 
\notag
\\
& = &
 \prod^{0\ \to\ 6}_i (\
\partial  - \ln'  ( \
\frac{ y_{7-i}  \prod_{s=1}^{6-i} T_s }{ y_{6-i} }) ) ,
\eean
The operator $D$ does not depend on the choice of
$\bs y\in P$. The kernel of $D$ is the $7$-dimensional space of
polynomials $V$ corresponding to $P$.

Let $\La_s$ be the unique dominant integral $A_6$-weight such that
$(\La_s,\al_i)=m_s^{(i)}$. Given a $6$-tuple of polynomials $\bs y$, define
the corresponding $A_6$-weight at infinity 
\be
\La_\infty=\sum_{s=1}^n\La_s - \sum_{i=1}^6(\deg y_i)\al_i.
\ee
The set of weights of infinity corresponding to elements of a
population $P$ form an orbit of Weyl group with respect to the shifted
action. 

Let $\bar\La_\infty$ be the unique dominant weight in the orbit of
$\La_\infty$ under the shifted action of the Weyl group.
It is conjectured in \cite{MV} that the number of $A_6$-populations
associated to $T_i$ for generic $z_j$ equals the multiplicity of
$L_{\bar\La_\infty}$ in $L_{\La_1}\otimes\dots\otimes L_{\La_n}$, where
  $L_\La$ is the irreducible $sl_7$-module of highest weight $\La$.

The population $P$ is isomorphic to the space of full flags in $V$,
the isomorphism is given by $F\to \bs y$ where
$y_i=W^\dagger(F_i)$. 
See Section 5 of \cite{MV} for details and proofs.

\subsection{$C_3$-populations}
Consider the root system of type $C_3$. Let $\al_3$ be the long root
and $\al_1, \al_2$ the short ones. We have 
\bea
(\al_1,\al_1)=2,\qquad (\al_2,\al_2)=2, \qquad (\al_1,\al_2)=-1,\\
(\al_2,\al_3)=-2, \qquad (\al_1,\al_3)=0, \qquad (\al_3,\al_3)=4.
\eea
The root system $C_3$ corresponds to the symplectic Lie algebra $sp_6$.
The Weyl group of $sp_6$ is generated by simple reflections
$s_i$, $i=1,2,3$. The shifted action of Weyl group on $C_3$-weights is
\be
s_i\cdot \la=\la-\langle\la+\rho,\al_i^\vee\rangle \al_i,
\ee 
where $\langle\rho,\al_i^\vee\rangle=1$ for all $i$.

Fix polynomials $T_i(x)=\prod_{s=1}^n(x-z_s)^{m_s^{(i)}}$ where
$i=1,2,3$, $x$ is a formal variable and $z_i,m_s^{(i)}$ are parameters.

A triple of polynomials $\bs y=(y_1,y_2,y_3)$ is called {\it
  generic} if $y_i$ have no multiple roots and if $y_i,y_{i+1}$ have no
  common roots. A generic triple of polynomials $\bs y$ is called
  {\it a Bethe triple} of  type $C_3$ 
 if there exist polynomials $\tilde
  y_1,\tilde y_2,\tilde y_3$ such that 
\be
W(y_1,\tilde y_1)=T_1y_2,\qquad W(y_2,\tilde y_2)=T_2y_1y_3^2,\qquad
  W(y_3,\tilde y_3)=T_3y_2.
\ee

If $\bs y$ is a Bethe triple and the triple $\bs
y^{(i)}=(y_1,\dots,\tilde y_i, \dots,y_3)$ is generic then $y^{(i)}$
is also a Bethe triple of $C_3$ type, see Theorem 3.7 in \cite{MV}.
We call $\bs y^{(i)}$ {\it an immediate descendent of $\bs y$ in the
direction $i$}.

Let $\bar P$ be the minimal set of Bethe triples which contains
$\bs y$ and such that for all $\tilde{ \bs y}\in\bar P$ all immediate
descendants  $\tilde {\bs y}^{(i)}$ are also in $\bar P$. 

The degrees of all coordinates of all triples of
$\bar P$ are simultaneously bounded.
{\it The $C_3$-population of Bethe triples originated at $\bs y$} 
is the closure of $\bar P$ in $(\Pb (\C[x]))^3$.

\bigskip

The set of population of Bethe triples associated with polynomials
$T_i$ is in one-to-one correspondence with the set of $7$-dimensional
self-dual spaces of polynomials with no base points and ramification 
polynomials $T_1,T_2,T_3,T_3,T_2,T_1$.
Given a population $P$ and $\bs y\in P$, we define a $6$-tuple $\bs
y^A$ by 
\be
\bs y^A=(y_1,y_2,y_3^2,y_3^2,y_2,y_1).
\ee
Then the self-dual space of polynomials $V$ which corresponds to $P$
is the kernel of $D(\bs y^A)$, given by \Ref{Diff}.

Let $\La_s$ be the unique dominant integral $C_3$-weight such that
$\langle\La_s,\al_i^\vee\rangle=m_s^{(i)}$. 
Given a triple of polynomials $\bs y$,
define the 
corresponding $C_3$-weight at infinity by
\be
\La_\infty=\sum_{s=1}^n\La_s - \sum_{i=1}^3(\deg y_i)\al_i.
\ee
The set of weights of infinity corresponding to elements of a
population $P$ form an orbit of Weyl group with respect to the shifted
action. 

Let $\bar\La_\infty$ be the unique dominant weight in the orbit of
$\La_\infty$ under the shifted action of the Weyl group.
It is conjectured in \cite{MV} that the number of $C_3$-populations
associated to $T_i$ for generic $z_j$ equals the multiplicity of
$L_{\bar\La_\infty}$ in $L_{\La_1}\otimes\dots\otimes L_{\La_n}$, where
  $L_\La$ is the irreducible $sp_6$-module of highest weight $\La$.

The population $P$ is isomorphic to the space of isotropic flags in $V$,
the isomorphism is given by $F\to \bs y$ where
$y_i=W^\dagger(F_i),$ $i=1,2,3$. 
See Sections 6 and 7 of \cite{MV} for details and proofs.

\subsection{$G_2$-populations}
Consider the root system of type $G_2$. Let $\al_1$ be the long root
and $\al_2$ the short one. We have 
\be
(\al_1,\al_1)=6,\qquad (\al_1,\al_2)=-3, \qquad (\al_2,\al_2)=2.
\ee
The Weyl group of $G_2$ is generated by simple reflections
$s_1,s_2$. The shifted action of Weyl group on $G_2$-weights is
\be
s_i\cdot \la=\la-\langle\la+\rho,\al_i^\vee\rangle\al_i,
\ee 
where $\langle \rho,\al_i^\vee\rangle=1$ for all $i$.

Fix polynomials
$T_i(x)=\prod_{s=1}^n(x-z_s)^{m_s^{(i)}}$, where
$i=1,2$, $x$ is a formal variable and $z_i,m_s^{(i)}$ are parameters.

A pair of polynomials $\bs y=(y_1,y_2)$ is called {\it
  generic} if $y_i$ have no multiple roots and no common roots.
  A generic pair of polynomials $\bs y$ is called
  {\it a Bethe pair} of  type $G_2$ 
 if there exist polynomials $\tilde y_1,\tilde y_2$ such that 
\be
W(y_1,\tilde y_1)=T_1y_2,\qquad W(y_2,\tilde y_2)=T_2y_1^3.
\ee

If $\bs y$ is a Bethe pair and the pair $\bs
y^{(1)}=(\tilde y_1,y_2)$  (resp. $\bs
y^{(2)}=(y_1,\tilde y_2)$)  is generic then $\bs y^{(1)}$ (resp. $\bs y^{(2)}$) 
is also a Bethe pair of $G_2$ type, see Theorem 3.7 in \cite{MV}.
We call $\bs y^{(i)}$ {\it an immediate descendent of $\bs y$ in the
direction $i$}. 

Let $\bar P$ be the minimal set of Bethe triples which contains
$\bs y$ and such that for all $\tilde{ \bs y}\in\bar P$ all immediate
descendants  $\tilde {\bs y}^{(i)}$ are also in $\bar P$. 

The degrees of all coordinates of all tuples of
$\bar P$ are simultaneously bounded.
{\it The $G_2$-population of Bethe pairs originated at $\bs y$} 
is the closure of $\bar P$ in $(\Pb(\C[x]))^2$.

Let $\La_s$ be the unique dominant integral $G_2$-weight such that
$\langle\La_s,\al_i^\vee\rangle=m_s^{(i)}$. 
Given a pair of polynomials $\bs y$, define
the corresponding $G_2$-weight at infinity
\be
\La_\infty=\sum_{s=1}^n\La_s - \sum_{i=1}^2(\deg y_i)\al_i.
\ee
The set of weights of infinity corresponding to elements of a
population $P$ form an orbit of Weyl group with respect to the shifted
action, see \cite{MV}, Theorem 3.12.

Let $\bar\La_\infty$ be the unique dominant weight in the orbit of
$\La_\infty$ under the shifted action of the Weyl group.
It is conjectured in \cite{MV} that the number of $G_2$-populations
associated to $T_i$ for generic $z_j$ equals the multiplicity of
$L_{\bar\La_\infty}$ in $L_{\La_1}\otimes\dots\otimes L_{\La_n}$, where
  $L_\La$ is the irreducible $G_2$ module of highest weight $\La$.

\bigskip

Given a $G_2$-population $P$ and $\bs y\in P$, we define a 3-tuple $\bs
y^C$ by 
\be
\bs y^C=(y_1,y_2,y_1).
\ee
The following lemma is obvious but very useful.

\begin{lem}
A pair of polynomials $\bs y$ is a Bethe pair of type $G_2$ associated
to polynomials $T_1,T_2$ if and
only if the triple of polynomials $\bs y^C$ is a Bethe triple of type
$C_3$ associated to polynomials $T_1,T_2,T_1$. Moreover, we have an
inclusion of populations $P(\bs y)\to P(\bs y^C)$ mapping $\tilde{\bs
  y}\mapsto \tilde{\bs y}^C$. $\Box$
\end{lem}

Assign to a $G_2$-population $P$ the space $V=\on{Ker}(D(((\bs y)^C)^A))$,
where $\bs y\in P$ and 
$D$ is defined by \Ref{Diff}. Then $V$ is a self-dual 
7-dimensional space of polynomials which does not depend on the choice
of $\bs y\in P$. Moreover for any $\bs y\in P$, $y_1\in V$.

See \cite{MV} for details and proofs.

\subsection{The space of $G_2$-population}
We show that the space of the $G_2$-population is self-self-dual. 

Let $P$ be a $G_2$-population and $V$ the corresponding self-dual
space of dimension 7.
We start with a description of degrees of polynomials in $V$ first.
The degrees of coordinates of
pairs of polynomials in $P$ are in one-to-one correspondence with the
Weyl group orbit (with respect to the shifted action).
It follows that there exists a pair $\bs y\in P$
such that $\deg \tilde y_1>\deg y_1$ and $\deg \tilde y_2>\deg
y_2$. Such a pair corresponds to the dominant weight in the Weyl group
orbit. Denote 
\be
\deg y_1=a, \qquad \deg y_2=b, \qquad \deg T_1=t_1, \qquad \deg
T_2=t_2.
\ee
So we have $\deg (\bs y)=(a,b)$, $\deg(\bs
y_1^{(2)})=(a, 3a-b+1+t_2)$ and $\deg(\bs y^{(1)})=(b-a+1+t_1,b)$.
In particular our choice of $\bs y$ means
\bean\label{dom}
b+1+t_1>2a,\qquad 3a+1+t_2>2b.
\eean  

In what follows we consider successive descendents of
a $G_2$-pair. Our notation 
$\bs y^{(1)(2)(1)}$ mean
$((\bs y^{(1)})^{(2)})^{(1)}$, etc. 

\begin{lem}\label{degrees}
The degrees of polynomials in $V$ are: 
$a$, $-a+b+1+t_1$, $2a-b+2+t_1+t_2$, $3+2t_1+t_2$,
$-2a+b+4+3t_1+t_2$, $a-b+5+3t_1+2t_2$, $-a+6+4t_1+2t_2$.
\end{lem}
\begin{proof}
These are degrees of the first coordinates of the descendents:
$\bs y$, $\bs y^{(1)}$, $\bs y^{(2)(1)}$, $\bs
y^{(1)(2)(1)'}$, $\bs y^{(1)(2)(1)}$, $\bs
y^{(2)(1)(2)(1)}$, $\bs
y^{(1)(2)(1)(2)(1)}$. Here by $\bs
y^{(1)(2)(1)'}$ we mean a difference between two different
descendents. 

The only non-trivial part here is to find the degree of 
$y_1^{(1)(2)(1)'}$. Let $y_1^{(1)\pm}$ be
two different immediate descendants such that
$y_1^{(1)+}=y_1^{(1)-}+\alpha y_1$.
Similarly, we will write $(\ )^\pm$ for two different descendants
and denote the difference by $(\ )'$.
We have the following Wronskians, with equalities up to
signs.
$$
W(y_1,y_1^{(1)\pm})=T_1y_2,~
W(y_2^{(1)(2)\pm}, y_2) = T_2 (y_1^{(1)^\pm})^3,~
W(y_1^{(1)(2)(1)\pm},y_1^{(1)\pm})= T_1y_2^{(1)(2)\pm}.
$$
From the last equation we conclude
that 
$$
F=W(y_1^{(1)(2)(1)'}, y_1^{(1)+}) = T_1y_2^{(1)(2)'} - 
\alpha W(y_1^{(1)(2)(1)-},y_1).
$$
We compute the degree of $F$.
We calculate
$$
W(F, T_1y_2)=
T_1^2T_2 ((y_1^{(1)+})^3-(y_1^{(1)-})^3) -
\alpha W(W(y_1^{(1)(2)(1)-},y_1),W(y_1^{(1)-},y_1))
$$
$$=\alpha^3 T_1^2 T_2y_1^3 + 3\alpha^2 T_1^2T_2 y_1^2 y_1^{(1)-}
+3\alpha T_1^2T_2 y_1 (y_1^{(1)-})^2
$$$$
\pm (y_1)(y_1^{(1)-})^{-1} \alpha  W (W(y_1^{(1)(2)(1)-},y_1^{(1)-}),
W(y_1^{(1)-},y_1))
$$
$$
= \alpha^3 T_1^2 T_2y_1^3 + 3\alpha^2 T_1^2T_2 y_1^2 y_1^{(1)-}
+(3\pm 1)\alpha T_1^2 T_2y_1 (y_1^{(1)-})^2.
$$
Here we have used the Wronskian identity
$W(W(u_1,u_2),W(u_1,u_3))=W(u_1,u_2,u_3)u_1$ which holds for all
$u_1,u_2,u_3$.

Using \Ref{dom} we see that the last term has the highest
degree which means that the degree of 
$W(F,T_1y_2)$ equals $2t_1+t_2+a+2\deg y_1^{(1)}$.
This implies that the degree of $F$ equals 
$1-t_1-b +2t_1+t_2+a+2\deg y_1^{(1)}$ which in turn
implies that 
$$
\deg(y_1^{(1)(2)(1)'})= 1-\deg y_1^{(1)} + \deg F =
3+2t_1+t_2.
$$

Since all the degrees in the statement of the lemma 
are present and different (written
in the increasing order) by \Ref{dom}, the lemma is proved.

\end{proof}

\begin{lem}\label{6 dim}
Any $G_2$-population $P$ is an algebraic variety of dimension $6$. The
set of the first coordinates of pairs in $P$ coincides with the set of
isotropic vectors in $V$. 
\end{lem}
\begin{proof}
Population $P$ is an algebraic variety by Corollary 3.13 in \cite{MV}.

Denote $J=\{y_1\ |\ (y_1,y_2)\in P\}$. Denote  as before $I\subset V$ the set
 of isotropic vectors. Since the first coordinates of
pairs in a $G_2$-population are also the first coordinates of triples in a
$C_3$-population we have
 $J\subseteq I$.  In particular $\dim J\leq \dim I=5$.

Consider the obvious projection of
algebraic varieties $P\to J$. 
It is well known 
that for any polynomial $q(x)\in \C[x]$ the set of planes of
polynomials $H\in Gr(2,\C[x])$ such that the Wronskian $W(H)$ is a scalar
multiple of $q(x)$, is finite. E.g. this follows from
Theorem 3.18 in \cite{MV}. Applied to $q(x)=T_2y_1^3$, it shows
that the fibers 
of $P\to J$ are at most 1-dimensional and we have $\dim P\leq \dim
J+1=6$. 

Now we show the opposite inequality. 
We have the following chain of descendents and their degrees:
\bea
\deg(\bs y^{(1)})&=&(b-a+1+t_1,b),\\
\deg (\bs y^{(1)(2)})&=&(b-a+1+t_1,2b-3a+4+3t_1+t_2), \\
\deg(\bs y^{(1)(2)(1)})&=&(
b-2a+4+3t_1+t_2,2b-3a+4+3t_1+t_2),\\ 
\deg (\bs y^{(1)(2)(1)(2)})&=&( b-2a+4+3t_1+t_2,b-3a+9+6t_1+3t_2), \\
\deg(\bs y^{(1)(2)(1)(2)(1)})&=&(-a+6+4t_1+2t_2,b-3a+9+6t_1+3t_2),\\
\deg (\bs
y^{(1)(2)(1)(2)(1)(2)})&=&(-a+6+4t_1+2t_2,-b+10+6t_1+4t_2).
\eea
 
We prove that the dimension of
the set of the descendents in the above chain increases by one. 
We recall that all our pairs are considered as elements in $\Pb\C[x]^2$.

The point is that given a pair $(y_1,y_2)$ on line $k$ we 
can uniquely recover the pair on the line $(k-1)$ that 
gave rise to $(y_1,y_2)$. Indeed, one of the coordinates is
always unchanged, and the other is uniquely determined
by the  degree restriction. Indeed, they are obtained by doing 
the reverse move from $(y_1,y_2)$ which results in a smaller
degree than before. Such move is always at best unique.

On the other hand, for every $(y_1,y_2)$ on line $k$ there
are $\C$ ways to extend it to line $(k+1)$. This shows
that the dimension of the set of descendants on line $(k+1)$ 
is one bigger than that on line $k$.

Our chain has length 6 and therefore $\dim P\geq 6$. It follows  that
$\dim P=6$ and $\dim J=5$. Since $\Pb I$ is an irreducible algebraic
variety of dimension $5$ and $J\subset \Pb I$, we see that 
 $\bar J=\Pb I$. On the other hand $J$ is closed.
\end{proof}

\begin{theorem}\label{population is selfself}
The space $V$ is self-self-dual.
\end{theorem}
\begin{proof}
To check the condition of Definition \ref{selfself},
we consider
the sequence of $C_3$ reproductions of the triple 
$(y_1,y_2,y_1)$ in the directions
1,2,3. Namely we have polynomials $\tilde y_i$ such that
\be
W(y_1,\tilde y_1)=T_1y_2,\quad  W(y_2,\tilde y_2)=T_2\tilde
y_1y_1^2,\quad W(\tilde y_3,y_1)=T_1\tilde y_2.
\ee
Then the triple $(\tilde y_1,\tilde y_2,\tilde y_3)$ is in a 
$C_3$-population. 
And in particular $\tilde y_3$ is a divided Wronskian of an
isotropic 3-space. 

\medskip

To show $\tilde y_3\in V$ we act on $\tilde y_3$ by the differential
operator $D(y_1,y_2,y_3^2,y_3^2,y_2,y_1)$. Applying the right factors
$\tilde y_3$ becomes (up to a constant) $\tilde y_2/y_1$ then
$y_1\tilde y_1/y_2$, then $1$ and then $0$. Therefore $\tilde y_3$ is
in $V$ which by definition is the kernel of this differential operator.

\medskip

Next we show that a generic $v\in V$ is $\tilde y_3$ for a suitable
$G_2$-pair $(y_1, y_2)$. Since everything is algebraic we just have to
show that the set of such $v$ (up to a scalar multiple) has dimension
6. We have 6-dimensional $G_2$-population. There exists 
a six-dimensional subset of the 
variety of $G_2$ pairs, such that each 
produces a variety of $\tilde y_3$ which has dimension at least
3. Indeed, the dimension is exactly 3 if we start with the dominant
pair $(y_1,y_2)$ satisfying $\deg y_1=a$, $\deg y_2=b$ as in Lemma
\ref{degrees}. To show that the dimension is at least three for 
a generic $G_2$ pair of the type ${\bf y}^{(1)(2)(1)(2)(1)(2)}$
of Lemma \ref{6 dim},
observe that at each step of the reproduction we can keep the 
descendant arbitrarily close to the previous pair, in terms of 
the projective space, by adding a big multiple of $y_i$. Hence,
there is an open set in the space of  ${\bf y}^{(1)(2)(1)(2)(1)(2)}$
which is contained in a neighborhood of $(y_1,y_2)$. 
Since the dimension 3 condition is clearly an open one, we see
that every $G_2$ pair which is close enough to $(y_1,y_2)$ satisfies it.

Now fix $\tilde y_3$. We need to show that the dimension of the variety
of $G_2$-pairs $(y_1,y_2)$ which produce $\tilde y_3$ is at most 3. Note, that
without loss of generality we can assume that  $\tilde y_3$ is
non-isotropic. Indeed, if we start with the dominant element of the
$G_2$-population $(y_1,y_2)$ then 
$\deg \tilde y_3=3+2t_1+t_2$ and $\tilde y_3$ is non-isotropic.
Since the condition of being is non-isotropic is open, it holds 
for the general descendants of the $G_2$ pairs in the neighborhood
of $(y_1,y_2)$.

Let $F$ be the
isotropic flag which corresponds to the $C_3$-triple $(y_1,y_2,y_1)$
and let $u_1,\dots,u_7$ be a Witt basis corresponding to $F$, such 
that  we have 
$u_1=y_1$, $W^\dagger(u_1,u_2)=y_2$ and $W^\dagger(u_1,u_2,u_3)=y_1^2$.

It follows that there exist constants $C_1,C_2,C_3$ such that the
$C_3$-triple $(\tilde y_1,\tilde y_2,\tilde y_3)$ corresponds to the
flag $\tilde F$, which is related to the basis
\be
\{\tilde u_1=u_1+C_1u_2,\tilde u_2=u_2+C_2u_3,\tilde
u_3=u_3+C_3u_4+C_3^2u_5/2,\tilde u_4,\tilde u_5,\tilde u_6,\tilde u_7\}
\ee
see Lemmas 6.14-6.16 in \cite{MV}. In particular, $W^\dagger(\tilde
u_1,\tilde u_2,\tilde u_3)=\tilde y_3^2$. 

There are only finitely many 3-spaces of polynomials with Wronskian
$\tilde y_3^2T_1^2T_2$. The 3-space $U$ spanned by $\tilde u_1,\tilde
u_2,\tilde u_3$ is
one of them. Note that the intersection of $F_3$ and $U$ contains $\tilde
u_1,\tilde u_2$ and therefore is at least 2-dimensional.

Now we observe that for any 2-dimensional subspace $U_1\subset U$ 
there are at most finitely many isotropic 3-spaces $F_3$ such that
$U_1\subset F_3\cap U$ and $\sqrt{W^\dagger(F_3)}\in F_3$. Indeed, the 
family of isotropic 3-spaces $F_3$ such that $U_1\subset F_3\cap U$ is
isomorphic to a
non-degenerate conic in a $\Pb (U_1^\perp/U_1)\cong \Pb^2$. 
Note that the space  $U$ belongs to this family and does not
contain $\sqrt{W^\dagger(U)}=\tilde y_3$ because $y_3$ is
non-isotropic. Therefore our family is a proper subvariety in a
non-degenerate conic, thus it is a finite set of points.

The dimension of the variety of 2-planes in $U$ is 2 and each $F_3$
produces a family of $G_2$-pairs of dimension 1, since 
$y_1=\sqrt{W^\dagger(F_3)}$ up to a scalar.
It follows that the dimension of the variety of $G_2$-pairs which
produce a given $\tilde y_3$ 
is at most 3. Thus the variety of all possible $\tilde y_3$ has
dimension at least 6 as needed.

\medskip

Finally, we show that  $\tilde y_3$ is orthogonal to $U$.
Recall that $\tilde u_1=\tilde y_1$,  $W(\tilde u_1,\tilde
u_2)=T_1\tilde y_2$. 
We denote $d=W(\tilde u_1,\tilde u_3)/T_1$, then
$W(\tilde y_2,d)=T_2\tilde y_1\tilde y_3^2$.

We calculate (up to constants):
$$
W(\tilde u_1,\tilde  u_2,\tilde u_3,\tilde y_3)=\tilde y_1^{-2}
W(W(\tilde y_1,\tilde u_2), W(\tilde y_1,\tilde u_3), 
W(\tilde y_1,\tilde y_3))
$$
$$
=
\tilde y_1^{-2}W(T_1\tilde y_2,T_1d,W(\tilde y_1,\tilde y_3))
=
T_1^3\tilde y_1^{-2}W(\tilde y_2,d,T_1^{-1}W(\tilde y_1,\tilde y_3))
$$
$$
=
T_1^3\tilde y_1^{-2}\tilde y_2^{-1}W(W(\tilde y_2,d),
W(\tilde y_2,T_1^{-1}W(\tilde y_1,\tilde y_3)))
$$
$$
=
T_1^3\tilde y_1^{-2}\tilde y_2^{-1}W(T_2\tilde y_1\tilde y_3^2,
T_1^{-2}W(W(y_1,\tilde y_3),W(\tilde y_1,\tilde y_3)))
$$
$$
=
T_1^3\tilde y_1^{-2}\tilde y_2^{-1}W(T_2\tilde y_1\tilde y_3^2,
T_1^{-2}\tilde y_3y_1^{-1}W(W(y_1,\tilde y_3),W(y_1,\tilde y_1)))
$$
$$
=
T_1^3\tilde y_1^{-2}\tilde y_2^{-1}W(T_2\tilde y_1\tilde y_3^2,
T_1^{-2}\tilde y_3y_1^{-1}W(T_1\tilde y_2,T_1y_2))
$$
$$
=
T_1^3\tilde y_1^{-2}\tilde y_2^{-1}W(T_2\tilde y_1\tilde y_3^2,
\tilde y_3y_1^{-1}T_2\tilde y_1y_1^2)
=T_1^3T_2^2\tilde y_2^{-1}\tilde y_3^2W(\tilde y_3,
y_1)
$$
$$
=T_1^4T_2^2 \tilde y_3^2=T_1^4T_2^2 W^\dagger(\tilde u_1,\tilde
u_2,\tilde u_2).
$$
The orthogonal complement $U^\perp$ of $U$ is the unique 
4-dimensional subspace of $V$ such that $U\subset U^\perp$ and
$W^\dagger (U^\perp)=W^\dagger(U)$. Therefore $U^\perp$ is the 
span of $u_1,u_2,u_3,\tilde y_3$. In particular $\tilde y_3\perp U$.
This shows that $V$ satisfies the condition of 
Definition \ref{selfself}.
\end{proof}

Theorem \ref{population is selfself}
shows that $V$ has the skew-symmetric 3-form by Theorem \ref{wrform}.
Thus for each element $\bs y$ of population $P$ we have a 
$G_2$-isotropic flag
$F$ in $V$ such that 
\be
F_1=W^\dagger(F_6)=y_1,\qquad
W^\dagger(F_2)=W^\dagger(F_5)=y_2,
\qquad W^\dagger(F_3)=W^\dagger(F_4)=y_1^2.
\ee

\begin{theorem}
The population $P$ is isomorphic to the variety $F^{\prp}(V)$ 
of $G_2$-isotropic flags in $V$. In particular $V$ is isomorphic to
the flag variety of group $G_2$. 
\end{theorem}
\begin{proof}
We already know that $P\subset F^{\prp}(V)$. In addition $\dim
P=\dim F^{\prp}(V)=6$. The variety 
$F^{\prp}(V)$ is a $\Pb_1$-bundle over an irreducible conic of
isotropic vectors and therefore is irreducible.  It follows that
$P=F^{\prp}(V)$.

The space $F^{\prp}(V)$ is isomorphic the $G_2$-isotropic flag variety by
Proposition \ref{G/B}.
\end{proof}
Let $P$ be a $G_2$-population and $\bs y\in P$ such that $\La_\infty$
is dominant integral. Let $F^\infty$ be the unique full flag of
$V$ such that the degrees of polynomials in $F_i$  are not larger that the
degrees of polynomials in $F_{i+1}$.
\begin{prop}
The closure of all elements of $P$ with the weights at
infinity equal $w\cdot \La_\infty$ is the Bruhat cell $G_w^{F^\infty}$.
\end{prop}
\begin{proof}
Completely parallel to the proof of Corollary 5.23 in \cite{MV}.
\end{proof}

\begin{theorem}
The set of $G_2$-populations associated to polynomials $T_1,T_2$ 
is in one to one correspondence with the set of self-self-dual spaces
of polynomials of dimension $7$ with ramification polynomials
$T_1,T_2,T_1,T_1,T_2,T_1$.
\end{theorem}
\begin{proof}
Consider a $G_2$-isotropic flag $F$. Then if the pair
$(F_1,W^\dagger(F_2))$ is generic then it is a 
$G_2$ Bethe pair.
Therefore, we only have to show that each self-self-dual space contains a
$G_2$-isotropic flag $F$ such that $(F_1,W^\dagger(F_2))$
form a generic pair. It is parallel to the proof of Theorems 7.5 and
7.10 in \cite{MV}.
\end{proof}

\bigskip

\subsection{Another description of self-self-dual spaces}
We show that Definition \ref{selfself} is equivalent to a simple
condition on 3-Wronskians.
\begin{theorem}\label{alter}
A self-dual space of polynomials is self-self-dual if and only if
\be
\{u^2 \ |\ u\in V\}=\{W^\dagger(u_1,u_2,u_3)\ |\ u_i\in V,\ (u_i,u_j)=0\}.
\ee
\end{theorem}
\begin{proof}
The only if part is Corollary \ref{2 to 1 wr}.

We need to show the if part. Let $V$ be a self-dual space. Recall that
such a space corresponds to a $C_3$-population.
Let $u_1,\dots,u_7$ be a basis  of $V$
such that $\deg u_i<\deg u_{i+1}$ and
$B(u_i,u_j)=(-1)^{i+1}\delta_i^{8-j}$. Then 
$W^\dagger (u_1,u_2,u_3)$ has the smallest degree among all divided
3-Wronskians and therefore we get $W^\dagger
(u_1,u_2,u_3)=cu_1^2$. It follows that $(y_1,y_2)$, where $y_1=u_1$, 
$y_2=W^\dagger(u_1,u_2)$ has the
reproduction properties of a $G_2$-pair. Namely there exist
 $\tilde
y_1=u_2$ and $\tilde y_2=c_1W^\dagger(u_1,u_3)$ such that $W(y_1,\tilde
y_1)=y_2T_1$ and $W(y_2,\tilde y_2)=y_1^3T_2$. The 
pair $(y_1,y_2)$ may be not generic and therefore it is not a
$G_2$-pair in general.

However the triples
$(y_1,\tilde y_2, y_1)$ and $(\tilde y_1, y_2,\tilde y_1)$ are clearly
in the same $C_3$-population and therefore correspond to some
isotropic flags. In
particular, these two triples again have the reproduction properties
of a $G_2$-pair. 

It follows that the condition of Definition \ref{selfself} is
satisfied by the same argument as in Theorem \ref{population is selfself}.
\end{proof}

\medskip 

\subsection{Examples of self-self-dual spaces}\label{monex}
The simplest example of a self-self-dual space of
polynomials is the space of polynomials $V$ of degree at most $6$. 
This space clearly corresponds to the population originated at the $G_2$-pair
$(1,1)$ where $T_1=T_2=1$. 

More generally, for every pair of integers $m<n$ a space spanned by monomials
$1,x^m,x^n,x^{m+n},x^{2m+n},x^{2n+m},x^{2m+2n}$
is self-self-dual. It corresponds to the population which 
originates at $(1,1)$ where $T_1=x^{m-1}$ and $T_2=x^{n-m-1}$.

\section {Standard bases of self-self-dual
spaces}\label{bases}

We recall that every self-dual space $V$ has a basis 
$\{v_1,\dots,v_7\}$  such that the divided 6-Wronskians are 
explicitly given by 
\be
W^\dagger(v_1,\dots,\hat v_i,\dots,v_7)=v_{8-i}.
\ee
Such a basis is called a Witt basis.
The basis $\{v_i\}$ is a Witt basis if and only if
$B(v_i,v_j)=(-1)^{i+1}\delta_i^{8-j}$, see \cite{MV}.

In this section we show that every self-self-dual space
has a Witt basis $\{v_1,\dots,v_7\}$ such that the 
all divided 3-Wronskians  are given explicitly as
a certain explicit linear combination of $v_iv_j$, see Table 1 below.
We call such a basis a standard basis. A Witt basis is a standard
basis if and only if the 3-form has the standard form, see Corollary
\ref{st 3-form} below.

\subsection{Standard dominant bases}
We keep the notation of Lemma \ref{degrees}. We also introduce 
$$m=-2a+b+1+t_1,~n=a-b+2+t_1+t_2.$$
Then $m<n$ and the degree list of Lemma
\ref{degrees} can be written as 
$$
a,a+m,a+n,a+m+n,a+2m+n,a+m+2n,a+2m+2n.
$$

First, we need we choose a Witt
basis of $V$ which is compatible with degrees. 
\begin{lem}\label{expl.witt}
The space $V$ has a Witt basis $\{v_1,\ldots, v_7\}$ 
of degrees 
$$
a,\ a+m,\ a+n,\ a+m+n,\ a+2m+n,\ a+m+2n,\ a+2m+2n
$$
and leading coefficients
\bea
&&1,\quad \frac 1m,\quad \frac 1{n(n-m)},\quad \frac 1{(m+n)nm},\quad
\frac 1{(2m+n)(m+n)(2m)m},\\
&&\frac 1{(m+2n)(2n)(m+n)n(n-m)},\quad
\frac 1{(2m+2n)(m+2n)(2m+n)(m+n)mn}
\eea
respectively.
\end{lem}

\begin{proof}
By Lemma 6.6, Lemma 6.7 in \cite{MV}, there exists a Witt
basis with the above degrees. Then we scale the basis
elements so that the leading coefficients are as above.
It is easy to check that the leading terms of the 6-Wronskians
are equal to the leading terms of the corresponding basis elements and
therefore after the scaling we again obtain a Witt basis.
\end{proof}

A Witt basis $\{v_1,\dots,v_7\}$ of $V$ with the above 
degrees and leading coefficients is called a {\it standard dominant
  basis} if 
\be
W^\dagger(v_1,v_5,v_6)=\frac 14 v_4^2,\qquad
W^\dagger(v_2,v_3,v_7)=\frac 12 v_4^2.
\ee

\begin{lem}\label{standard.basis}
Let $V$ be self-self-dual. Then there exists a standard dominant basis of $V$ 
with degrees and leading  terms specified by Lemma \ref{expl.witt}.
\end{lem}

\begin{proof}
By Corollary \ref{2 to 1 wr}, there are two ways of 
writing $v_4^2$ as a divided Wronskian of an isotropic 
three-space, up to a constant. We call the corresponding
spaces $U_1$ and $U_2$. For each $U_i$ we can 
find a basis $\{f_{i,1},f_{i,2},f_{i,3}\}$ of increasing degrees. 
Moreover, without loss of generality 
we can assume that the leading coefficients and the degrees of the 
$f_i$ are among those of Lemma \ref{expl.witt}. 
Since $f_{i,j}$ are isotropic their degrees
are not equal to the degree of $v_4$.

It is easy to see that we must have 
$$
\deg f_{i,1}+\deg f_{i,2}+\deg f_{i,3} = 3a +3m +3n.
$$
This implies that for each $i$ the degrees are either 
$(a,a+2m+n,a+m+2n)$ or $(a+m,a+n,a+2m+2n)$. 
It is impossible to have degrees from the first list for both $i$
since that would imply that $U_1\cap U_2\ni v_1$. 
Similarly, if both $U_i$ had degrees from the second list,
we would have $\dim U_i\cap {\rm Span}(v_1,v_2,v_3)=2$,
which again implies $U_1\cap  U_2\neq 0$.
As a result, we can assume that 
$f_{1,1},f_{2,1},f_{2,2},v_4,f_{1,2},f_{1,3},f_{2,3}$ have
degrees and leading terms of Lemma \ref{expl.witt}.

The above basis is a a Witt basis provided that 
$B(f_{1,i},f_{2,j})=0$ if $i+j\neq 3$ and 
$B(f_{1,1},f_{2,3})=B(f_{1,2},f_{2,2})=-B(f_{1,3},f_{2.1})=1$.

Most of these equalities hold automatically. For example, 
$f_{1,2}=v_5+...$ and $f_{2,2}=v_3+...$
where the dots denote linear combinations of 
$v_i$ with lower $i$. This implies that
$B(f_{1,2},f_{2,2})=1$ in view of the pairing of $v_i$.

There are exactly three equalities that are not true for a generic
choice of $f_{i,j}$. Namely, we may not have
$$
B(f_{1,2},f_{2,3})=0, \quad B(f_{1,3},f_{2,2})=0, \quad B(f_{1,3},f_{2,3})=0.
$$
Then we  get $B(f_{1,3},f_{2,2})=0$
by adding an appropriate scalar multiple of $f_{2,1}$ to $f_{2,2}$.
Note that such an addition operation does not change the leading term.

Similarly, to get $B(f_{1,2},f_{2,3})=0$ and $B(f_{1,3},f_{2,3})=0$ we
add an appropriate linear combination of $f_{2,1}$
and $f_{2,2}$ to $f_{2,3}$. This finishes the proof.
\end{proof}

\begin{example}
Let $V$ be the space of polynomials of degree at most 6. Recall that $V$ is
self-self-dual. Then the basis
\be
\{1,\ x,\ \frac{x^2}{2!},\ \frac{x^3}{3!},\ \frac{x^4}{4!},\
\frac{x^5}{5!},\ \frac{x^6}{6!}\}
\ee
is standard dominant.
\end{example}

\subsection{3-Wronskians in standard bases}
Fix a standard dominant basis $\{v_1,\dots,v_7\}$ in $V$. 
We identify the spinor space $\hat V$ with the 
polynomials in odd variables $\hat v_5$, $\hat v_6$ and
$\hat v_7$ as before. We recall that 
$v_5,v_6,v_7$ act by multiplications by $\hat v_5,\hat
v_6, \hat v_7$ respectively, $v_1$, $v_2$, $v_3$ act
by differentiations $-\partial/\partial \hat v_7$, $\partial/\partial
\hat v_6$, $-\partial/\partial \hat v_5$ 
respectively, and $v_4$ acts by $\frac 1{\sqrt 2}(-1)^{\rm deg}$,
where ${\rm deg}$ is the degree of the odd polynomial.

Formula \eqref{explicitrho} allows us
to fix the explicit form of the pairing $\hat B$.
Namely,
$$
\hat B(\hat 1, \hat v_5\hat v_6
\hat v_7)
=\hat B(\hat v_6,\hat v_5\hat v_7)=1,
\quad \hat B(\hat v_7,\hat v_5\hat v_6)= 
\hat B(\hat v_6\hat v_7,\hat v_5)=-1,
$$
and all other pairings of basis elements are zero.

Recall that by Proposition \ref{there is a p} there is 
an invariant surjection $\psi: \hat V\to V$ such that
if $\psi(\rho(U))=f$ then the divided Wronskian of $U$  
is proportional to $f^2$.

\begin{lem}\label{expl.p}
The point $p$ corresponding to the invariant surjection 
$\psi$ is given by the formula 
$p=\hat v_5\hat v_6 +\frac1{\sqrt{2}} \hat v_7$.
\end{lem}
\begin{proof}
The polynomial $v_4^2$ is proportional to the divided Wronskians
of the spans of $(v_1,v_5,v_6)$ and $(v_2,v_3,v_7)$.
Applying the spinor embedding $\rho$ we obtain the lines
$\C \hat v_5\hat v_6$ and $\C \hat v_7$ in $\hat V$.
Since the invariant surjection is a projection from a point $p$, 
we see that, up to a constant multiple, 
$$
p=\hat v_5\hat v_6 +\alpha \hat v_7
$$
for some $\alpha \in \C$. 

Consider the isotropic vector $v_1+\beta v_2$.
Its square can be uniquely up to constant written as a
divided Wronskian of the annihilator of $(v_1+\beta v_2)\cdot p$.
Since 
$$(v_1+ \beta v_2)\cdot p =
(-\partial/\partial \hat v_7 +\beta \partial/\partial \hat v_6)
\cdot (\hat v_5\hat v_6 +\alpha \hat v_7) = -\alpha\hat 1 - \beta \hat v_5,$$
the corresponding isotropic 3-space is spanned
by $v_1,v_2, \alpha^2 v_3 +\sqrt 2 \alpha\beta v_4 +\beta^2v_5$.
Therefore, we must have 
$$
W^\dagger (v_1,v_2,\alpha^2 v_3 +\sqrt 2 \alpha\beta v_4 
+\beta^2v_5) =  c(\beta) (v_1+\beta v_2)^2.
$$
Comparing the leading coefficients of both sides, we obtain
that $c(\beta)$ is a constant. Moreover, we must have 
$$
\alpha^2 W^\dagger(v_1,v_2,v_3) = cv_1^2,\quad
\alpha W^\dagger(v_1,v_2,v_4) = \sqrt 2c v_1v_2,\quad
W^\dagger(v_1,v_2,v_5) = c v_2^2
$$
for some constant $c$.
Comparing the leading 
coefficients on the both sides of the last two equations, we obtain
the equalities
\be
\alpha \frac {m(m+n)n}{m(m+n)mn}= c{\sqrt 2} \frac 1 m, \qquad
\frac {m(2m+n)(m+n)}{m(2m+n)(m+n)(2m)m}=c\frac1{m^2},
\ee
which give $c=\frac 12$, $\alpha = \frac 1{\sqrt 2}$.
\end{proof}

The explicit knowledge of $p$ and $\hat B$
allows us to calculate all divided 3-Wronskians in the standard
dominant basis.
The next theorem is the main result of this section.

\begin{theorem}\label{all35}
The Wronskians of the basis elements of a standard dominant basis
of $V$ are given in the Table 1 below. 
\end{theorem}

\begin{proof}
The theorem essentially amounts to the calculation of the map $\phi$
of Proposition \ref{phinu}, since by 
Corollary \ref{general wrons} the divided Wronskian map to $\C[x]$ is 
a composition of the map $\frac12\phi$ and the 
multiplication map $Sym^2 V\to \C[x]$.

Since we have fixed $p=\hat v_5\hat v_6 +\frac 1{\sqrt 2}\hat v_7$
and $\hat B$, by Proposition \ref{phinu} $\phi= C\nu$.
From the definition of the standard dominant basis, we have 
$\phi(v_1\wedge v_2\wedge v_3) = v_1\otimes v_1$.
We compute $\nu(v_1\wedge v_2 \wedge v_3)$:
\bea
\lefteqn{\langle \nu^*(v_1\wedge v_2\wedge
v_3),(\sum_{i=1}^7\alpha_iv_i)\otimes (\sum_{i=1}^7\alpha_iv_i)\rangle=}\\
&&=\hat B\left(v_1v_2v_3(\sum_{i=1}^7\alpha_i v_i)\cdot
 (\hat v_5\hat v_6 +\frac 1{\sqrt 2}\hat v_7),\
(\sum_{i=1}^7\alpha_i v_i) \cdot  (\hat v_5\hat v_6 +
\frac 1{\sqrt 2}\hat v_7)\right) \\
&&=\hat B\left(\alpha_7\hat 1, 
(\sum_{i=1}^7\alpha_i v_i) \cdot  (\hat v_5\hat v_6 +
\frac 1{\sqrt 2}\hat v_7)\right)
=\hat B(\alpha_7\hat 1,\alpha_7\hat v_5\hat v_6\hat v_7)=
\alpha_7^2=(B(v_1,\sum_{i=1}^7\alpha_i v_i))^2
\eea
Therefore $\nu(v_1,v_2,v_3) = v_1\otimes v_1$, so $C=1$
and $\phi=\nu$.

It is now routine to calculate all divided Wronskians. We
calculate $W^\dagger(v_1,v_4,v_7)$ as an example 
and leave the rest to the reader.
\bea
\lefteqn{
\langle \nu^*(v_1\wedge v_4\wedge v_7),\
(\sum_{i=1}^7\alpha_iv_i)\otimes(\sum_{i=1}^7\alpha_iv_i)\rangle=}\\
&& =\hat B\left(v_1v_4v_7(\sum_{i=1}^7\alpha_i v_i)\cdot
 (\hat v_5\hat v_6 +
 \frac 1{\sqrt 2}\hat v_7),
(\sum_{i=1}^7\alpha_i v_i) \cdot  (\hat v_5\hat v_6 +
\frac 1{\sqrt 2}\hat v_7)\right)\\
&&=\hat B(-\alpha_1 \frac 12\hat 1+\alpha_2\frac 1{\sqrt 2}\hat v_5
+\alpha_3\frac 1{\sqrt 2}\hat v_6+\alpha_4\frac 1{2}
\hat v_5\hat v_6, \\
&&-\alpha_1\frac 1{\sqrt 2} \hat 1 
- \alpha_2 \hat v_5 
- \alpha_3 \hat v_6 +\alpha_4(\frac 1{\sqrt 2}\hat v_5\hat v_6
-\frac 12\hat v_7) 
+\alpha_5 \frac 1{\sqrt 2} \hat v_5\hat v_7 
+\alpha_6 \frac 1{\sqrt 2} \hat v_6\hat v_7
+\alpha_7\hat v_5\hat v_6\hat v_7)\\
&&
=-\frac 12\alpha_1\alpha_7 - \frac 12\alpha_2\alpha_6
+\frac 12\alpha_3\alpha_5+\frac 14\alpha_4^2.
\eea
This implies that 
$$\phi(v_1\wedge v_4\wedge v_7)=
\frac 14(-v_1\otimes v_7 -v_7\otimes v_1 
-v_2\otimes v_6 -v_6\otimes v_2 + v_3\otimes v_5
+v_5\otimes v_3 + v_4\otimes v_4).$$
Hence $W^\dagger (v_1,v_4,v_7)=
-\frac 12v_1v_7-\frac 12v_2v_6+\frac 12 v_3v_5+\frac 14v_4^2$.
\end{proof}

\begin{figure}
\noindent
\begin{tabular}{|c|c|c|c|c|}
\hline
123 & 124 & 125& 126& 127
\\ 
$v_1^2$
& $v_1v_2$
& $\frac 12 v_2^2$
& $-\frac 12 v_1v_4 +\frac 12 v_2v_3$
& $-v_1v_5 + \frac 12 v_2v_4$
\smallskip
\\ \hline
134&135&136&137&145
\\ 
$v_1v_3$
&
$ \frac 12 v_1v_4 +\frac 12 v_2v_3$
&
$\frac 12 v_3^2$
&
$-v_1v_6 + \frac 12 v_3v_4$
&
$\frac 12 v_2v_4$
\smallskip
\\ \hline 
146&147&156&157&167
\\ 
$\frac 12 v_3v_4$
&
$
\begin{array}{c}{-\frac 12v_1v_7-\frac 12v_2v_6}\\{+\frac 12 v_3v_5+\frac 14v_4^2}\end{array}
$
&
$ \frac 14 v_4^2$
&
$ -\frac 12 v_2v_7 +\frac 12 v_4v_5$
&
$-\frac 12 v_3v_7 +\frac 12 v_4v_6$

\smallskip
\\ \hline
234&235&236&237&245
\\ 
$v_1v_4$
&
$v_1v_5+\frac 12 v_2v_4$
&
\;\;$v_1v_6+\frac 12 v_3v_4$\;\;
&
$\frac 12 v_4^2$
&
$ v_2v_5$
\smallskip \\ \hline 
246&247&256&257&267
\\ 
$
\noindent\begin{array}{c}
 {\frac 12 v_1v_7+\frac 12 v_2v_6}\\
 {+\frac 12 v_3v_5+\frac 14 v_4^2}
\end{array}$
&
$v_4v_5$
&
$\frac 12 v_2v_7+\frac 12 v_4v_5$
&
$v_5^2$
&
$-\frac 12 v_4v_7+v_5v_6$
\smallskip
\\   \hline
345&346&347&356&357 
\\
$
\noindent\begin{array}{c}{-\frac 12 v_1v_7+\frac 12 v_2v_6
}\\
{+\frac 12 v_3v_5-\frac 14 v_4^2}
\end{array}$
&
$ v_3v_6$
&
$v_4v_6$
&
$\frac 12 v_3v_7+\frac 12 v_4v_6$
&
$\frac 12 v_4v_7 + v_5v_6$
\smallskip
\\  \hline
367&456&457&467&567
\\ 
$v_6^2$
&
$\frac 12 v_4v_7$
&
$ v_5v_7$
&
$v_6v_7$
&
$\frac 12 v_7^2$
\smallskip
\\  \hline
\end{tabular}

\medskip

\centerline{
{\bf Table 1.} The entry under $ijk$ is $W^\dagger(v_i,v_j,v_k)$.}
\end{figure}


A basis $\{v_1,\dots,v_7\}$ 
of $V$ is called {\it standard} if the divided 3-Wronskians
$W^\dagger(v_i,v_j,v_k)$ are given in Table 1. A standard dominant
basis is standard.

\begin{prop}\label{66}
Let $V$ be a self-dual space with a standard basis. Then $V$
is self-self-dual. 
\end{prop}
\begin{proof}
Define $\phi_V:\ \La^3V\to Sym^2 V$ by the formulas in Table 1. 
Let $U$ be any self-self-dual space with a standard basis
$\{u_1,\dots,u_7\}$ and the map $\phi_U:\ \La^3U\to Sym^2 U$. 
Define the map $\iota:\ U\to V$ sending
$u_i\mapsto v_i$. Then it induces the maps $\La^3U\to \La^3 V$ and
$Sym^2U\to Sym^2 V$ which obviously intertwine $\phi_U$ and $\phi_V$. 
The image of any isotropic 3-space in $U$ is a tensor
square of an element in $U$. Therefore 
the image of any isotropic 3-space in $V$ is a tensor
square of an element in $V$. It follows that divided Wronskian of 
any isotropic 3-space in $V$ is a a square of an element of $V$. Now
the proposition follows from Theorem \ref{alter}.
\end{proof}

For a self-self-dual space $V$ the knowledge of only a few $3$-Wronskians
is sufficient to decide whether a Witt basis $\{v_1,\ldots,v_7\}$
is self-dual.
\begin{lem}\label{smallcheck}
A Witt basis $\{v_1,\ldots,v_7\}$ of a self-self-dual
space is standard if and only if
$$
W^\dagger(v_2,v_3,v_7)=\frac 12 v_4^2,~
W^\dagger(v_1,v_5,v_6)=\frac 14 v_4^2,~
$$
$$W^\dagger(v_1,v_2,v_3)=v_1^2,~
W^\dagger(v_1,v_2,v_4)=v_1v_2,~
W^\dagger(v_1,v_2,v_5)=\frac 12 v_2^2.
$$
\end{lem}

\begin{proof}
The if part is a tautology. To prove the only if part,
notice that the first two Wronskians assure that $p=\hat v_5\hat v_6+
\alpha v_7$. As in the proof of Lemma \ref{expl.p} 
we look at 
$$
W^\dagger (v_1,v_2,\alpha^2 v_3 +\sqrt 2 \alpha\beta v_4 
+\beta^2v_5) =  c(\beta) (v_1+\beta v_2)^2.
$$
Since $v_1^2,v_1v_2,v_2^2$ are linearly independent,
the last three Wronskians assure $\alpha=\frac 1{\sqrt 2}$.
Then the argument of Theorem \ref{all35} goes through, since
it only uses $W^\dagger(v_1,v_2,v_3)=v_1^2$ to fix the constant.
\end{proof}

\subsection{3-form in a standard basis}
Now we read off the explicit formula for the trilinear form $w$
of Theorem \ref{wrform} in the standard basis $\{v_i\}$.
\begin{prop}\label{expl 3-form}
For $i<j<k$ we have $w(v_i,v_j,v_k)=0$ except for 
\bea
w(v_2\wedge v_4\wedge v_6)=w(v_1\wedge v_4\wedge v_7)=
-w(v_3\wedge v_4\wedge v_5) = \\
=-w(v_1\wedge v_5\wedge v_6)= -{\frac{1}{2}} \ w(v_2\wedge v_3\wedge v_7)=
\frac14\ .
\eea
\end{prop}
\begin{proof}
Since $w$ and $w_\psi$ are 
$G_2(p)$-invariant, by Lemma \ref{only inv}
$w$ is a constant multiple of $w_\psi$.
Therefore, 
there exists a constant $C$, such that
\be
w(a\wedge b\wedge c)= 
C\hat B(abc\cdot p,p).
\ee

As we saw in the proof of Lemma \ref{smallcheck},
$p=\hat v_5\hat v_6 +\frac 1{\sqrt 2} \hat v_7$, so
we need to calculate 
$$\frac 1{\sqrt 2} (\hat B(abc\cdot \hat v_5\hat v_6,\hat v_7)
+\hat B(abc\cdot\hat v_7,\hat v_5\hat v_6))
+\hat B(abc\cdot \hat v_5\hat v_6,\hat v_5\hat v_6)
+\frac 12\hat B(abc\cdot \hat v_7,\hat v_7).$$

Since $\hat B(q_1,q_2)$ is zero unless the degrees of $q_i$
add up to three,
the first of these terms is nonzero only if the degree
of $abc$ is zero, which can only happen if it is of the
form $v_{\leq 3}v_{4}v_{\geq 5}$. It is also easy to see that 
the term is zero unless $abc$ is of the form $v_iv_4v_{8-i}$,
for which it equals $(-\frac 12)$ for $i=1,2$ and $\frac 12$
for $i=3$.

The last two terms are nonzero only for $abc=v_2v_3v_7$
and $abc=v_1v_5v_6$ respectively, when they equal 
$1$ and $\frac 12$ respectively.  Finally, from 
$W^\dagger(v_1,v_5,v_6)=\frac 14 v_4^2$ and
$B(v_4,v_4)=-1$ we obtain the constant: $C=-\frac 12$.
\end{proof}

Note that the value $w(a\wedge b\wedge c)$ can be computed by applying
$B$ to the corresponding element in Table 1. For example,
\be
w(v_1\wedge v_4 \wedge v_7)=-\frac12 B(v_1,v_7)-\frac12 B(v_2,v_6)+
\frac12 B(v_3,v_5)+\frac14 B(v_4,v_4)=\frac14.
\ee


\begin{cor}
Let $w_1$ be any non-degenerate skew-symmetric 3-form in a 7-dimensional
space $V$ associated to a non-degenerate bilinear form $B$. 
Then there is a basis $\{v_1,\dots,v_7\}$ of $V$ such that
$B(v_i,v_j)=(-1)^{i+1}\delta_i^{8-j}$ and 
$w_1$ is a scalar
multiple of the form described in Proposition \ref{expl 3-form}.
\end{cor}
\begin{proof}
By Proposition 
\ref{stab=G2}, the group $SL(V)$ acts transitively on the set of
all non-degenerate forms considered up to a constant. That implies the
corollary. 
\end{proof}

\medskip

Next we describe the set of standard bases in a self-self-dual
space $V$.

\begin{prop}\label{action on bases}
The group $G_2(p)$ acts transitively on the set of standard bases.
\end{prop}
\begin{proof}
The group $G_2(p)$ acts on the set of standard bases by Proposition
\ref{g2v}. Any two standard bases can be mapped to each other by an
element $g$ of the orthogonal group $SO(V)$. But then $g$ preserves the
3-form and therefore belongs to $G_2(p)$ by Lemma
\ref{st=G2}. Therefore the action is transitive.
\end{proof}

\begin{cor}\label{st 3-form}
A basis $\{v_1,\dots, v_7\}$ of a self-self-dual space
 $V$ is standard if and only if
$B(v_i,v_j)=(-1)^{i+1}\delta_{i}^{8-j}$ and 
the 3-form has the form
described in Proposition \ref{expl 3-form}.
\end{cor}
\begin{proof}
Follows from Proposition \ref{action on bases}.
\end{proof}

\begin{cor}\label{basis-flag} Let basis $\{v_1,\dots, v_7\}$ in $V$ be
  standard and let
  $F_i$ be the span of $v_1,\dots,v_i$. Then $F$ is a $G_2$-isotropic flag.
\end{cor}
\begin{proof}
Follows from Proposition \ref{action on bases}.
\end{proof}

\begin{cor}
For every $G_2$-isotropic flag $F\in \mc F^{\prp}$, there exists 
a standard basis 
$\{v_1,\dots,v_7\}$, such that $v_1,\dots,v_i$ span $F_i$. 
\end{cor}
\begin{proof}
The group $G_2$ acts transitively on $\mc F^{\prp}$ by Proposition
\ref{G/B} and on the set of standard bases by Proposition \ref{action
  on bases}. One standard basis exists by Lemma
\ref{standard.basis}. It is related to a $G_2$-isotropic flag by Corollary
\ref{basis-flag}. Then this basis is mapped by
$G_2$-action to a standard basis related to any other given
$G_2$-isotropic flag.
\end{proof}

We finish with a description of  the reproduction procedure in a standard
basis.

Let $\{v_1,\dots,v_7\}$ be a standard basis and let $(y_1,y_2)$ be the
 corresponding element of the $G_2$-population, $y_1=v_1$, 
$y_2=W^\dagger(v_1,v_2)$.
\begin{prop}
For any $c\in\C$ the two bases of $V$ given by 
\begin{align*} 
&\{v_1+cv_2,v_2,v_3+2cv_4+2c^2v_5,v_4+2cv_5,v_5,v_6+cv_7,v_7\},\\
&\{v_1,v_2+cv_3,v_3,v_4,v_5+cv_6,v_6,v_7\}
\end{align*}
are standard. The corresponding 
flags are in a bijective correspondence with
the set of the immediate descendents $(\tilde y_1,y_2)$ 
and $(y_1,\tilde y_2)$ 
of $(y_1,y_2)$ in the first and in the second directions respectively.
\end{prop}
\begin{proof}
It is straightforward to check that the
new bases are standard using Table 1 and Lemma \ref{smallcheck}.
To check that they correspond to the descendants, the only non-trivial 
check is $W(y_2,\tilde y_2)=T_2 y_1^3$ (up to a constant)
for the second direction. Since 
$y_2= W^\dagger(v_1,v_2)$ and 
$\tilde y_2=W^\dagger(v_1,v_2+cv_3)$,
we have 
$$
W(y_2,\tilde y_2)=T_1^{-2}W(W(v_1,v_2),W(v_1,v_2+cv_3))=
T_1^{-2}v_1W(v_1,v_2,v_2+cv_3)
$$$$
=cT_1^{-2}T_1^2T_2v_1^3
=cT_2v_1^3
=cT_2y_1^3.
$$
\end{proof}

\end{document}